\documentclass{article}
\usepackage{amssymb}

\usepackage{amsmath}

\newtheorem{theorem}{Theorem}

\newtheorem{corollary}[theorem]{Corollary}

\newtheorem{definition}[theorem]{Definition}
\newtheorem{example}[theorem]{Example}

\newtheorem{lemma}[theorem]{Lemma}

\newtheorem{proposition}[theorem]{Proposition}
\newtheorem{remark}[theorem]{Remark}

\newenvironment{proof}[1][Proof]{\textbf{#1.} }{\ \rule{0.5em}{0.5em}}

\begin{document}

\title{Spaces of Infinite Measure and Pointwise Convergence of the Bilinear
Hilbert and Ergodic Averages Defined by $L^{p}$-Isometries}
\date{September 1, 2007; revised January 31, 2008.
The second author was supported by NSF Grant DMS-0556389.}
\author{Earl Berkson \\
Department of Mathematics\\
University of Illinois\\
1409 W. Green St\\
Urbana, IL 61801 U.S.A\\
berkson@math.uiuc.edu \and Ciprian Demeter \\
School of Mathematics\\
Institute for Advanced Study\\
Einstein Drive\\
Princeton, NJ 08540 U.S.A\\
demeter@math.ias.edu}
\maketitle

\begin{abstract}
We generalize the respective ``double recurrence'' results of Bourgain and
of the second author, which established for pairs of $L^{\infty }$ functions
on a finite measure space the a.e. convergence of the discrete bilinear
ergodic averages and of the discrete bilinear Hilbert averages defined by
invertible measure-preserving point transformations. Our generalizations are
set in the context of arbitrary sigma-finite measure spaces and take the
form of a.e. convergence of such discrete averages, as well as of their
continuous variable counterparts, when these averages are defined by
Lebesgue space isometries and act on $L^{p_{1}}\times L^{p_{2}}$ ($%
1<p_{1},p_{2}<\infty $, $p_{1}^{-1}+p_{2}^{-1}<3/2$). In the setting of an
arbitrary measure space, this yields the a.e.\thinspace\ convergence of
these discrete bilinear averages when they act on $L^{p_{1}}\times L^{p_{2}}$
and are defined by an invertible measure-preserving point transformation.
\end{abstract}

\section{Introduction\label{sec1}}

For an arbitrary measure space $\left( X,\sigma \right) $, we shall denote
by $\mathcal{A}\left( \sigma \right) $ the algebra under pointwise
operations consisting of all complex-valued $\sigma $-measurable functions
on $X$ (identified modulo equality $\sigma $-a.e. on $X$). The class of all
real-valued functions belonging to $\mathcal{A}\left( \sigma \right) $ will
be denoted by $\Re \left( \mathcal{A}\left( \sigma \right) \right) $. In
this setting we shall use the following terminology and notation.

\begin{definition}
\label{Define_Averages}Let $T$ be a linear bijection of $\mathcal{A}\left(
\sigma \right) $ onto $\mathcal{A}\left( \sigma \right) $. For each real
number $r\geq 1$, we shall denote the integer part of $r$ by $\left[ r\right]
$, and for $f\in \mathcal{A}\left( \sigma \right) $, $g\in \mathcal{A}\left(
\sigma \right) $, we define the corresponding discrete bilinear ergodic
average $A_{r,T}\left( f,g\right) $ and the corresponding discrete bilinear
Hilbert average $H_{r,T}\left( f,g\right) $ by writing pointwise on $X$,%
\begin{eqnarray}
A_{r,T}\left( f,g\right) &=&\frac{1}{\left[ r\right] }\sum_{n=0}^{\left[ r%
\right] -1}\left( T^{n}f\right) \left( T^{-n}g\right) \text{;}  \label{E1.1}
\\
H_{r,T}\left( f,g\right) &=&\sum_{0<\left| n\right| \leq \left[ r\right] }%
\frac{\left( T^{n}f\right) \left( T^{-n}g\right) }{n}\text{.}  \label{E1.2}
\end{eqnarray}
\end{definition}

\bigskip Our principal concern will be to generalize the double recurrence
theorem of Bourgain for discrete bilinear ergodic averages \cite{Bou} and
its counterpart for discrete bilinear Hilbert averages (recently established
in Theorem 1.2 of \cite{Dem}), whose statements are reproduced as the
following theorem.

\begin{theorem}
\bigskip \label{Jean_Ciprian}Suppose that $\left( \mathfrak{X}\text{,}\rho
\right) $ is a finite measure space, and $\phi $ is an invertible
measure-preserving point transformation of $\left( \mathfrak{X}\text{,}\rho
\right) $\ onto $\left( \mathfrak{X}\text{,}\rho \right) $. Let $\mathcal{T}%
: $ $\mathcal{A}\left( \rho \right) \rightarrow \mathcal{A}\left( \rho
\right) $ denote composition with $\phi $. Then for every $f\in L^{\infty
}\left( \rho \right) $, and every $g\in L^{\infty }\left( \rho \right) $,
each of the sequences $\left\{ A_{k,\mathcal{T}}\left( f,g\right) \right\}
_{k=1}^{\infty }$ and $\left\{ H_{k,\mathcal{T}}\left( f,g\right) \right\}
_{k=1}^{\infty }$ converges $\rho $-a.e. on $\mathfrak{X}$ to a
corresponding function belonging to $\mathcal{A}\left( \rho \right) $.
\end{theorem}

Our main result, which is stated as follows, generalizes Theorem \ref%
{Jean_Ciprian} in the direction of $L^{p}$-isometries for sigma-finite
measure spaces. (See Theorem \ref{one-parameter-ver} in \S \ref{sec4} below
for the continuous variable version of this generalization.)

\begin{theorem}
\label{Main_Result}Suppose that $\left( \Omega ,\mu \right) $ is a
sigma-finite measure space, and let $U$ be a bijective linear mapping of $%
\mathcal{A}\left( \mu \right) $ onto $\mathcal{A}\left( \mu \right) $ such
that the following two conditions hold.

\begin{enumerate}
\item[(i)] Whenever $\left\{ g_{k}\right\} _{k=1}^{\infty }\subseteq 
\mathcal{A}\left( \mu \right) $, $g\in \mathcal{A}\left( \mu \right) $, and $%
g_{k}\rightarrow g$ $\mu $-a.e. on $\Omega $, it follows that as $%
k\rightarrow \infty $, $U$ $\left( g_{k}\right) \rightarrow U$ $\left(
g\right) $ $\mu $-a.e. on $\Omega $, and $U^{-1}$ $\left( g_{k}\right)
\rightarrow U$ $^{-1}\left( g\right) $ $\mu $-a.e. on $\Omega $.

\item[(ii)] The restriction $U\left| L^{p}\left( \mu \right) \right. $ is a
surjective linear isometry of $L^{p}\left( \mu \right) $ onto $L^{p}\left(
\mu \right) $ for $0<p\leq \infty $.
\end{enumerate}

\noindent Suppose further that%
\begin{gather}
1<p_{1},p_{2}<\infty \text{;}  \label{E1.3} \\
\frac{1}{p_{1}}+\frac{1}{p_{2}}=\frac{1}{p_{3}}<\frac{3}{2}\text{.}
\label{E1.4}
\end{gather}

\noindent Then for every $f\in L^{p_{1}}\left( \mu \right) $, and every $%
g\in L^{p_{2}}\left( \mu \right) $, each of the sequences%
\begin{equation*}
\left\{ A_{k,U}\left( f,g\right) \right\} _{k=1}^{\infty }\text{, }\left\{
H_{k,U}\left( f,g\right) \right\} _{k=1}^{\infty }
\end{equation*}
converges $\mu $-a.e. on $\Omega $ and in the metric topology of $%
L^{p_{3}}\left( \mu \right) $ to a corresponding function belonging to $%
\mathcal{A}\left( \mu \right) $.
\end{theorem}

Although the proof of Theorem \ref{Main_Result} will be deferred to \S \ref%
{sec3}, the hypotheses of Theorem \ref{Main_Result} on $\left( \Omega ,\mu
\right) $ and $U$ will remain in effect henceforth. Theorem \ref{Main_Result}
has the following corollary, which is valid for arbitrary measure spaces,
and which likewise obviously implies Theorem \ref{Jean_Ciprian}, since in
the setting of any finite measure space $\left( \mathfrak{X}\text{,}\rho
\right) $, the inclusion $L^{\infty }\left( \rho \right) \subseteq
L^{2}\left( \rho \right) $ holds. (A variant of this result, likewise valid
for all measure spaces, is described below for the continuous variable
averages in Corollary \ref{1param4all_meas-Cor}.)

\begin{corollary}
\label{4all_measure_spaces}Suppose that $\left( X,\sigma \right) $ is an
arbitrary measure space, and let $\tau $ be an invertible measure-preserving
point transformation of $\left( X,\sigma \right) $ onto $\left( X,\sigma
\right) $. Suppose also that $p_{1}$, $p_{2}$, $p_{3}$ satisfy (\ref{E1.3})
and (\ref{E1.4}), and let $f\in L^{p_{1}}\left( \sigma \right) $, $g\in
L^{p_{2}}\left( \sigma \right) $. Then each of the sequences%
\begin{equation*}
\left\{ \dfrac{1}{k}\sum_{n=0}^{k-1}f\left( \tau ^{n}\right) g\left( \tau
^{-n}\right) \right\} _{k=1}^{\infty }\text{, }\left\{ \sum_{0<\left|
n\right| \leq k}\dfrac{f\left( \tau ^{n}\right) g\left( \tau ^{-n}\right) }{n%
}\right\} _{k=1}^{\infty }
\end{equation*}
converges $\sigma $-a.e. on $X$ and in the metric topology of $%
L^{p_{3}}\left( \sigma \right) $ to a corresponding function belonging to $%
\mathcal{A}\left( \sigma \right) $.
\end{corollary}

\begin{proof}
Since $f\in L^{p_{1}}\left( \sigma \right) $, we can write $\left\{ x\in
X:\left| f\left( x\right) \right| >0\right\} =\bigcup_{j=1}^{\infty }E_{j}$,
where $\sigma \left( E_{j}\right) <\infty $, for each $j\in \mathbb{N}$.
Putting%
\begin{equation*}
Y=\bigcup_{n\in \mathbb{Z}}\bigcup_{j=1}^{\infty }\tau ^{n}\left(
E_{j}\right) \text{,}
\end{equation*}%
we see that for all $n\in \mathbb{Z}$: $\tau ^{n}\left( Y\right) =Y$, and $%
f\left( \tau ^{n}\left( x\right) \right) =0$ for all $x\in X\setminus Y$. So
in order to establish the desired $\sigma $-a.e. convergence on $X$ it
suffices to prove that for $\sigma $-almost all $x\in Y$, each of the
sequences $\left\{ \dfrac{1}{k}\sum_{n=0}^{k-1}f\left( \tau ^{n}\left(
x\right) \right) g\left( \tau ^{-n}\left( x\right) \right) \right\}
_{k=1}^{\infty }$ and $\left\{ \sum_{0<\left| n\right| \leq k}\dfrac{f\left(
\tau ^{n}\left( x\right) \right) g\left( \tau ^{-n}\left( x\right) \right) }{%
n}\right\} _{k=1}^{\infty }$ is convergent. But this follows immediately
upon application of Theorem \ref{Main_Result} to the sigma-finite measure
space $\left( Y,\sigma \right) $ and the composition operator corresponding
to the restriction $\tau $ $\left| Y\right. $.

Likewise for convergence in $L^{p_{3}}\left( X\text{,}\sigma \right) $,
which can also be seen as follows. By combining the reasoning regarding $Y$
with the maximal estimates in the setting of sigma-finite measure spaces of
Theorems 9 and 10 in \cite{BBCG} (whose statements are reproduced in Theorem %
\ref{bisublinear_maxml} below), we see without difficulty that by dominated
convergence each of the sequences 
\begin{equation*}
\left\{ \dfrac{1}{k}\sum_{n=0}^{k-1}f\left( \tau ^{n}\right) g\left( \tau
^{-n}\right) \right\} _{k=1}^{\infty }\text{, }\left\{ \sum_{0<\left|
n\right| \leq k}\dfrac{f\left( \tau ^{n}\right) g\left( \tau ^{-n}\right) }{n%
}\right\} _{k=1}^{\infty }
\end{equation*}%
converges in the metric of $L^{p_{3}}\left( \sigma \right) $.
\end{proof}

\begin{remark}
\label{Direcns}While our main concern will be with extensions of Theorem \ref%
{Jean_Ciprian} to spaces of infinite measure, we briefly comment here on
aspects in which Theorem \ref{Main_Result} generalizes Theorem \ref%
{Jean_Ciprian}, when Theorem \ref{Main_Result} is restricted to the finite
measure space setting. If $\left( \Omega _{0},\mu _{0}\right) $ is a finite
measure space, then a bijective linear mapping $U$ of $\mathcal{A}\left( \mu
_{0}\right) $ onto $\mathcal{A}\left( \mu _{0}\right) $ satisfying
conditions (i) and (ii) of Theorem \ref{Main_Result} need not have any of
its odd powers implemented by an invertible measure-preserving point
transformation of $\Omega _{0}$ onto $\Omega _{0}$. Such an example is
furnished by \S 343J of \cite{Frem}, which, taken in conjunction with
standard considerations about measure-preserving set transformations (see
pgs. 452-454 of \cite{Doob}), furnishes a complete, non-atomic, finite
measure space $\left( \tilde{\Omega}_{0},\tilde{\mu}_{0}\right) $ \ whose
measure algebra is a separable metric space, together with a self-inverse
algebra automorphism $U_{0}$ of $\mathcal{A}\left( \tilde{\mu}_{0}\right) $
onto $\mathcal{A}\left( \tilde{\mu}_{0}\right) $ satisfying conditions (i)
and (ii) of Theorem \ref{Main_Result}, but is such that $U_{0}\left|
L^{\infty }\left( \tilde{\mu}_{0}\right) \right. $ cannot be expressed by
composition with an invertible measure-preserving point transformation of
the measure space $\left( \tilde{\Omega}_{0},\tilde{\mu}_{0}\right) $ onto $%
\left( \tilde{\Omega}_{0},\tilde{\mu}_{0}\right) $. So Theorem \ref%
{Jean_Ciprian} does not directly apply here, although in this particular
example, since $U_{0}$ is self-inverse, the $\tilde{\mu}_{0}$-a.e.
convergence from Theorem \ref{Main_Result} will hold trivially--in fact,
whenever $f\in \mathcal{A}\left( \tilde{\mu}_{0}\right) $ and $g\in \mathcal{%
A}\left( \tilde{\mu}_{0}\right) $. (For information about the general
relationships between measure-preserving set transformations and
measure-preserving point transformations of non-atomic finite measure spaces
satisfying separability conditions, see \S 41 of \cite{Halmos} and \cite%
{PRH-vN}.) In the special case of Theorem \ref{Main_Result} where the
setting is an arbitrary finite measure space $\left( \Omega _{0},\mu
_{0}\right) $ and $\left( f,g\right) \in L^{\infty }\left( \mu _{0}\right)
\,\times \,L^{\infty }\left( \mu _{0}\right) $, the convergence $\mu _{0}$%
-a.e. of the averages $\left\{ A_{k,U}\left( f,g\right) \right\}
_{k=1}^{\infty }$ and $\left\{ H_{k,U}\left( f,g\right) \right\}
_{k=1}^{\infty }$ can be deduced directly from Theorem \ref{Jean_Ciprian} in
conjunction with maximal results from \cite{BBCG} (which are quoted as
Theorem \ref{bisublinear_maxml} below) by following G.-C. Rota's
``dilation'' theory for measure spaces (\cite{Doob1},\cite{Rota}) in a
spirit similar to that in Chapter IV, \S 4, of \cite{Stein}. We omit the
details of this line of reasoning for the finite measure space setting,
since we will be proving Theorem \ref{Main_Result} in full generality by
following a different path.
\end{remark}

\bigskip The reasoning used in \cite{Dem} to deduce the a.e. convergence of
the discrete averages $H_{k,\mathcal{T}}\left( f,g\right) $ under the
conditions of Theorem \ref{Jean_Ciprian} above also furnished a new proof of
the result of \cite{Bou} \ for the a.e. convergence of the discrete averages 
$\left\{ A_{k,\mathcal{T}}\left( f,g\right) \right\} _{k=1}^{\infty }$ in
the same circumstances. Our strategy (particularly as regards infinite
measures) for treating the wider scope of Theorem \ref{Main_Result} and for
obtaining its continuous variable counterpart will be to combine suitable
modifications of the unified coverage of $\left\{ A_{k,\mathcal{T}}\left(
f,g\right) \right\} _{k=1}^{\infty }$ and $\left\{ H_{k,\mathcal{T}}\left(
f,g\right) \right\} _{k=1}^{\infty }$ in \cite{Dem} with the treatment of
the relevant bisublinear maximal operators in Theorems 9, 10, and 13 of \cite%
{BBCG}. Accordingly, the remaining four sections of this article will be
organized as follows. In \S \ref{sec2} we collect some background items that
furnish key structural tools for the demonstration in \S \ref{sec3} of
Theorem \ref{Main_Result}. In particular, \S \ref{sec2} \ includes an
expanded version of Lemma 3.1 of \cite{Dem} aimed at providing, in the form
of a suitable oscillation estimate, an abstract sufficiency criterion for $%
\mu $-a.e. convergence (see Lemma \ref{Demeter_Lemma3_1} \ below). The way
is then opened for proceeding from Lemma \ref{Demeter_Lemma3_1} and the
maximal bisublinear theorems of \cite{BBCG} to derive Theorem \ref%
{Main_Result} in \S \ref{sec3} after the development in \S \ref{sec2a} of
the relevant oscillation estimates, which will take the form of general
discretization and transference results. We close by treating the continuous
variable model in \S \ref{sec4}, which, in particular, establishes the
counterpart of Theorem \ref{Main_Result} for one-parameter groups of
Lebesgue space isometries associated with the arbitrary sigma-finite measure
space $\left( \Omega ,\mu \right) $ (see Theorem \ref{one-parameter-ver}
below).

\bigskip Henceforth, the following notation will be in effect. If $A$ is a
subset of a given set $Y$, then, except where otherwise indicated, the
characteristic function of $A$, defined on $Y$, will be designated by $\chi
_{A}$, and the restriction to $A$ of a function $F$\ defined on $Y$ will be
written $F\left| A\right. $. The collection of all mappings of a set $E$
into a set $W$ will be denoted by $W^{E}$. Lebesgue measure on $\mathbb{R}$
(respectively, counting measure on $\mathbb{Z}$) will be symbolized by $%
\mathfrak{m}_{\mathbb{R}}$ (respectively, by $\mathfrak{m}_{\mathbb{Z}}$).
Given an arbitrary measure space $\left( X,\sigma \right) $, and a function $%
f\in \mathcal{A}\left( \sigma \right) $, we shall denote by $\lambda \left(
f,\sigma ;\left( \cdot \right) \right) $ the \textit{distribution function }%
of $f$ specified by:%
\begin{equation}
\lambda \left( f,\sigma ;y\right) =\sigma \left( \left\{ x\in X:\left|
f\left( x\right) \right| >y\right\} \right) \text{, for each real number }y>0%
\text{,}  \label{E1.6}
\end{equation}%
and we shall follow standard notation by writing%
\begin{equation*}
\left\| f\right\| _{L^{1,\infty }\left( \sigma \right) }=\sup \left\{
y\lambda \left( f,\sigma ;y\right) :y\in \mathbb{R},y>0\right\} \text{.}
\end{equation*}%
Given a positive real number $\xi $ and a function $f:\mathbb{R\rightarrow C}
$, we shall symbolize by $\delta _{\xi }f$ the dilation of $f$ by $\xi $,
which is defined by%
\begin{equation}
\left( \delta _{\xi }f\right) \left( x\right) =\frac{1}{\xi }f\left( \frac{x%
}{\xi }\right) \text{, for all }x\in \mathbb{R}\text{.}  \label{e1.7}
\end{equation}%
The letter ``$C$'' with a (possibly empty) set of subscripts will signify a
constant which depends only on those subscripts, and which can change its
value from one occurrence to another.

\section{\protect\bigskip Background Items\label{sec2}}

The following proposition, a version of Corollary 3.1 in \cite{Lamp} (see
also Proposition 5 and Remark 5-(i) in \cite{BBCG}), describes the structure
of the operator $U$ in the hypotheses of Theorem \ref{Main_Result}.

\begin{proposition}
\label{structure_of_U}In the setting of the arbitrary sigma-finite measure
space $\left( \Omega ,\mu \right) $, let $U$ be a bijective linear mapping
of $\mathcal{A}\left( \mu \right) $ onto $\mathcal{A}\left( \mu \right) $
such that the conditions (i) and (ii) in the hypotheses of Theorem \ref%
{Main_Result} hold. Then there are unique sequences $\left\{ h_{j}\right\}
_{j=-\infty }^{\infty }$ and $\left\{ \Phi _{j}\right\} _{j=-\infty
}^{\infty }$ such that for each $j\in \mathbb{Z}$:

\begin{enumerate}
\item[(j)] $h_{j}\in \mathcal{A}\left( \mu \right) $, with $\left|
h_{j}\right| =1$ on $\Omega $, and $\Phi _{j}$ is an algebra automorphism of 
$\mathcal{A}\left( \mu \right) $ onto $\mathcal{A}\left( \mu \right) $;

\item[(jj)] for every $f\in \mathcal{A}\left( \mu \right) $, $U^{j}f$ is
expressed by the pointwise product on $\Omega $ of the functions $h_{j}$ and 
$\Phi _{j}\left( f\right) $;

\item[(jjj)] whenever $\left\{ f_{k}\right\} _{k=1}^{\infty }\subseteq 
\mathcal{A}\left( \mu \right) $, $f\in \mathcal{A}\left( \mu \right) $, and $%
f_{k}\rightarrow f$ $\mu $-a.e. on $\Omega $, it follows that as $
k\rightarrow \infty $, 
$$\Phi _{j}\left( f_{k}\right) \rightarrow \Phi
_{j}\left( f\right) $$ $\mu $-a.e. on $\Omega $.
\end{enumerate}

\hspace{0.5in}\hspace{0.5in}This unique sequence $\left\{ \Phi _{j}\right\}
_{j=-\infty }^{\infty }$ has the property that%
\begin{gather}
\mu \left( E\right) =\int_{\Omega }\,\Phi _{j}\left( \chi _{E}\right) \,d\mu 
\text{,}  \label{e2.1} \\
\hspace{0.5in}\hspace{0.5in}\hspace{0.5in}\hspace{0.5in}\text{for each }j\in 
\mathbb{Z}\text{, and each }\mu \text{-measurable set }E\text{.}  \notag
\end{gather}%
The sequences $\left\{ h_{j}\right\} _{j=-\infty }^{\infty }$ and $\left\{
\Phi _{j}\right\} _{j=-\infty }^{\infty }$ also enjoy the following
properties for arbitrary $j\in \mathbb{Z},$ $k\in \mathbb{Z}$.%
\begin{eqnarray}
\Phi _{j}\left( g\right) &\geq &0\text{, for each }g\geq 0\text{ belonging
to }\mathcal{A}\left( \mu \right) \text{.}  \label{e2.2} \\
\left| \Phi _{j}\left( f\right) \right| ^{\alpha } &=&\Phi _{j}\left( \left|
f\right| ^{\alpha }\right) \text{, \ for }f\in \mathcal{A}\left( \mu \right)
,\text{ and }0<\alpha <\infty \text{.}  \label{e2.3} \\
\Phi _{j+k}\left( f\right) &=&\Phi _{j}\left( \Phi _{k}\left( f\right)
\right) \text{, for every }f\in \mathcal{A}\left( \mu \right) .  \label{e2.4}
\\
h_{j+k}\left( x\right) &=&h_{j}\left( x\right) \left( \Phi _{j}h_{k}\right)
\left( x\right) \text{, for }\mu \text{-almost all }x\in \Omega .
\label{e2.5}
\end{eqnarray}
\end{proposition}

\begin{remark}
\label{Props_of_Phij}(a) It is clear that for each $j\in \mathbb{Z}$ the
restriction $\Phi _{j}\left| L^{p}\left( \mu \right) \right. $ is a
surjective linear isometry of $L^{p}\left( \mu \right) $ onto $L^{p}\left(
\mu \right) $, for $0<p\leq \infty $.

(b) By (\ref{e2.2}), $\Phi _{j}\left( g\right) $ is real-valued for each $%
j\in \mathbb{Z}$, and each real-valued $g$ belonging to $\mathcal{A}\left(
\mu \right) $. \ Moreover, from (\ref{e2.4}) we see that $\Phi _{-j}=\Phi
_{j}^{-1}$. Hence application of (\ref{e2.2}) to $\Phi _{j}$ and to $\Phi
_{-j}$ shows that the restriction of $\Phi _{j}$ to the class $\Re \left( 
\mathcal{A}\left( \mu \right) \right) $ consisting of the real-valued $\mu $%
-measurable functions on $\Omega $ is a surjective order isomorphism. It
follows that for each finite sequence $\left\{ g_{k}\right\}
_{k=1}^{N}\subseteq $ $\Re \left( \mathcal{A}\left( \mu \right) \right) $,%
\begin{equation}
\Phi _{j}\left( \sup_{1\leq k\leq N}\,g_{k}\right) =\sup_{1\leq k\leq
N}\,\Phi _{j}\left( \,g_{k}\right) \text{.}  \label{e2.5a}
\end{equation}

(c) It is readily seen from the foregoing considerations that each $\Phi
_{j} $ preserves distribution functions. That is, for $f\in \mathcal{A}%
\left( \mu \right) $, each positive real number $y$, and each $j\in \mathbb{Z%
}$, we have, in the notation of (\ref{E1.6}),%
\begin{equation}
\lambda \left( \Phi _{j}\left( f\right) ,\mu ;y\right) =\lambda \left( f,\mu
;y\right) \text{.}  \label{e2.6}
\end{equation}
\end{remark}

The following expanded version of the sufficiency criterion for a.e.
convergence in Lemma 3.1 of \cite{Dem} will play a pivotal role in our
considerations.

\begin{lemma}
\label{Demeter_Lemma3_1}Let $\left( \Omega ,\mu \right) $ be a sigma-finite
measure space, and suppose that $\left\{ f_{n}\right\} _{n=1}^{\infty }$ is
a sequence of complex-valued $\mu $-measurable functions defined on $\Omega $
which has the following property: there is a positive real constant $\Theta $
such that for every positive integer $J\geq 2$, and every sequence of
positive integers $u_{1}<u_{2}<\cdots <u_{J}$, we have 
\begin{equation}
\left\| \left\{ \sum_{j=1}^{J-1}\sup_{\substack{ n\in \mathbb{N},  \\ %
u_{j}\leq n<u_{j+1}}}\left| f_{n}-f_{u_{j+1}}\right| ^{2}\right\}
^{1/2}\right\| _{L^{1,\infty }\left( \mu \right) }\leq \Theta J^{1/4}\text{.}
\label{e2.7}
\end{equation}%
Then there is a $\mu $-measurable function $f:\Omega \rightarrow \mathbb{C}$
such that $\left\{ f_{n}\right\} _{n=1}^{\infty }$ converges to $f$ $\mu $%
-a.e. on $\Omega $.
\end{lemma}

\begin{proof}
In view of the sigma-finiteness of $\left( \Omega ,\mu \right) $ and the
form of the hypothesis (\ref{e2.7}), without loss of generality we can and
shall assume henceforth in the proof of this lemma that $\mu \left( \Omega
\right) <\infty $. Putting $L\left( x\right) =\limsup_{n\rightarrow \infty
}\left| f_{n}\left( x\right) \right| $, for all $x\in \Omega $, and writing $%
A=\left\{ x\in \Omega :L\left( x\right) =\infty \right\} ,$ we claim that $%
\mu \left( A\right) =0$. For each $j\in \mathbb{N}$, let%
\begin{equation*}
E_{j}=\left\{ x\in \Omega :\left| f_{1}\left( x\right) \right| \leq
j\right\} \text{.}
\end{equation*}%
For the claim that $\mu \left( A\right) =0$, it clearly suffices to show
that $\mu \left( A\bigcap E_{j}\right) =0$, for each $j\in \mathbb{N}$.
Assume to the contrary that $\mu \left( A\bigcap E_{j_{0}}\right) >0$, for
some $j_{0}\in \mathbb{N}$. Observe that by Egoroff's Theorem there is $%
B\subseteq A\bigcap E_{j_{0}}$ such that $\mu \left( B\right) >\dfrac{\mu
\left( A\bigcap E_{j_{0}}\right) }{2}$, and such that the sequence%
\begin{equation*}
\left\{ \min_{1\leq k\leq n}\exp \left( -\left| f_{k}\right| \right)
\right\} _{n=1}^{\infty }
\end{equation*}%
converges to the zero function uniformly on $B$. Accordingly, for each $m\in 
\mathbb{N}$, there is $N\in \mathbb{N}$ such that for all $n\geq N$,%
\begin{equation*}
\sup_{1\leq k\leq n}\left| f_{k}\right| >m+j_{0}\text{ on }B\text{.}
\end{equation*}%
For $1\leq \nu \leq N$, we have on $\Omega $, 
\begin{eqnarray}
\left| f_{\nu }-f_{1}\right| &\leq &\left| f_{\nu }-f_{N+1}\right| +\left|
f_{N+1}-f_{1}\right|  \label{e2.7I} \\
&\leq &2\sup_{1\leq k<N+1}\left| f_{k}-f_{N+1}\right| \text{,}  \notag
\end{eqnarray}%
and it now follows that%
\begin{equation*}
\mu \left( B\right) \leq \mu \left\{ x\in \Omega :\sup_{1\leq k<N+1}\left|
f_{k}-f_{1}\right| >m\right\} \leq \mu \left\{ x\in \Omega :\sup_{1\leq
k<N+1}\left| f_{k}-f_{N+1}\right| >m/2\right\} \text{.}
\end{equation*}%
This, together with an application of (\ref{e2.7}) (for $J=2$, $u_{1}=1$, $%
u_{2}=N+1$), shows that%
\begin{equation*}
\dfrac{m\mu \left( A\bigcap E_{j_{0}}\right) }{2}<m\mu \left( B\right) \leq
\Theta 2^{5/4}\text{, for all }m\in \mathbb{N}\text{.}
\end{equation*}%
Since $m\in \mathbb{N}$ is arbitrary, we can let $m\rightarrow \infty $ in
this to contradict the supposition that $\mu \left( A\bigcap
E_{j_{0}}\right) >0$, thereby establishing the claim that $\mu \left(
A\right) =0$. Hence%
\begin{equation*}
\mu \left( \Omega \setminus \bigcup_{k=1}^{\infty }\left\{ x\in \Omega
:\sup_{n\in \mathbb{N}}\left| f_{n}\left( x\right) \right| \leq k\right\}
\right) =0\text{,}
\end{equation*}%
and so by confining attention to each set $\left\{ x\in \Omega :\sup_{n\in 
\mathbb{N}}\left| f_{n}\left( x\right) \right| \leq k\right\} $ separately,
we can assume (in addition to (\ref{e2.7}) and $\mu \left( \Omega \right)
<\infty )$ that $\left\{ f_{n}\right\} _{n=1}^{\infty }\subseteq L^{\infty
}\left( \Omega ,\mu \right) $, with%
\begin{equation*}
\sup_{n\in \mathbb{N}}\left\| f_{n}\right\| _{L^{\infty }\left( \Omega ,\mu
\right) }<\infty \text{.}
\end{equation*}

From this point onwards, the proof of Lemma \ref{Demeter_Lemma3_1} is a
straightforward adaptation of the proof sketched in \cite{Dem} for Lemma 3.1
therein.
\end{proof}

\begin{remark}
\bigskip \label{re_lemmas_needed} Given a sequence $\left\{ f_{n}\right\}
_{n=1}^{\infty }\subseteq $ $\mathcal{A}\left( \mu \right) $, a positive
integer $J\geq 2$, and any sequence of positive integers $u_{1}<u_{2}<\cdots
<u_{J}$, we can adapt the elementary reasoning of (\ref{e2.7I}) to infer
readily that for $1\leq j\leq J-1$,%
\begin{equation*}
\sup_{\substack{ n\in \mathbb{N},  \\ u_{j}\leq n<u_{j+1}}}\left|
f_{n}-f_{u_{j+1}}\right| \leq 2\,\sup_{\substack{ n\in \mathbb{N},  \\ %
u_{j}<n\leq u_{j+1}}}\left| f_{n}-f_{u_{j}}\right| \text{.}
\end{equation*}%
Consequently,%
\begin{equation*}
\left\| \left\{ \sum_{j=1}^{J-1}\sup_{\substack{ n\in \mathbb{N},  \\ %
u_{j}\leq n<u_{j+1}}}\left| f_{n}-f_{u_{j+1}}\right| ^{2}\right\}
^{1/2}\right\| _{L^{1,\infty }\left( \mu \right) }\leq 2\,\left\| \left\{
\sum_{j=1}^{J-1}\sup_{\substack{ n\in \mathbb{N},  \\ u_{j}<n\leq u_{j+1}}}%
\left| f_{n}-f_{u_{j}}\right| ^{2}\right\} ^{1/2}\right\| _{L^{1,\infty
}\left( \mu \right) }\text{.}
\end{equation*}%
Similarly we also have%
\begin{equation}
\left\| \left\{ \sum_{j=1}^{J-1}\sup_{\substack{ n\in \mathbb{N},  \\ %
u_{j}<n\leq u_{j+1}}}\left| f_{n}-f_{u_{j}}\right| ^{2}\right\}
^{1/2}\right\| _{L^{1,\infty }\left( \mu \right) }\leq 2\,\left\| \left\{
\sum_{j=1}^{J-1}\sup_{\substack{ n\in \mathbb{N},  \\ u_{j}\leq n<u_{j+1}}}%
\left| f_{n}-f_{u_{j+1}}\right| ^{2}\right\} ^{1/2}\right\| _{L^{1,\infty
}\left( \mu \right) }\text{.}  \label{e2.8}
\end{equation}%
Hence the minorant in (\ref{e2.7}) can be equivalently replaced by the
minorant in (\ref{e2.8}).
\end{remark}

In order to cover the $\mu $-a.e. convergence for all $f\in L^{p_{1}}\left(
\mu \right) $ and all $g\in L^{p_{2}}\left( \mu \right) $ in the conclusion
of Theorem \ref{Main_Result}, we shall require the following bisublinear
maximal theorem for this setting (a combination of the statements of
Theorems 9 and 10 in \cite{BBCG}, which were obtained by discretizing and
transferring Michael Lacey's bisublinear maximal theorems for the classical
setting \cite{Lacey}).

\begin{theorem}
\label{bisublinear_maxml}Assume the hypotheses of Theorem \ref{Main_Result}
on $\left( \Omega ,\mu \right) $, $U$, $p_{1}$, $p_{2}$, and $p_{3}$. For $%
f\in L^{p_{1}}\left( \mu \right) $ and $g\in L^{p_{2}}\left( \mu \right) $,
define the maximal functions $\mathcal{H}_{U}\left( f,g\right) $ and $%
\mathcal{M}_{U}\left( f,g\right) $ by writing pointwise on $\Omega $, 
\begin{eqnarray}
\mathcal{H}_{U}\left( f,g\right) &=&\sup_{j\in \mathbb{N}}\left|
H_{j,U}\left( f,g\right) \right| =\sup_{j\in \mathbb{N}}\,\left|
\sum_{0<\left| n\right| \leq j}\frac{\left( U^{n}f\right) \,\left(
U^{-n}g\right) }{n}\right| \text{;}  \label{e2.12} \\
\mathcal{M}_{U}\left( f,g\right) &=&\sup_{j\in \mathbb{N}}\,\frac{1}{2j+1}%
\sum_{n=-j}^{j}\left| U^{n}f\right| \,\left| U^{-n}g\right| \text{.}
\label{e2.13}
\end{eqnarray}%
Then there are positive constants $\mathfrak{M}_{p_{1},p_{2}}$ and $%
\mathfrak{N}_{p_{1},p_{2}}$, each depending only on $p_{1}$ and $p_{2}$,
such that for all $f\in L^{p_{1}}\left( \mu \right) $ and all $g\in
L^{p_{2}}\left( \mu \right) $:%
\begin{eqnarray}
\left\| \mathcal{H}_{U}\left( f,g\right) \right\| _{L^{p_{3}}\left( \mu
\right) } &\leq &\mathfrak{M}_{p_{1},p_{2}}\,\left\| f\right\|
_{L^{p_{1}}\left( \mu \right) }\,\left\| g\right\| _{L^{p_{2}}\left( \mu
\right) }\text{;}  \label{e2.14} \\
\left\| \mathcal{M}_{U}\left( f,g\right) \right\| _{L^{p_{3}}\left( \mu
\right) } &\leq &\mathfrak{N}_{p_{1},p_{2}}\,\left\| f\right\|
_{L^{p_{1}}\left( \mu \right) }\,\left\| g\right\| _{L^{p_{2}}\left( \mu
\right) }\text{.}  \label{e2.15}
\end{eqnarray}
\end{theorem}

\section{\protect\bigskip Discretized and Transferred Oscillation Estimates%
\label{sec2a}}

In this section, we develop general discretization and transference results
that will implement the derivation of Theorem \ref{Main_Result} by setting
up suitable applications of Lemma \ref{Demeter_Lemma3_1} and Theorem \ref%
{bisublinear_maxml}. The first step in this process will be to discretize
the following oscillation theorem for the real line (Theorem 1.4 of \cite%
{Dem}), which plays a seminal role in our considerations.

\begin{theorem}
\label{Dem_Thm_1_4}Let $K:$ $\mathbb{R\rightarrow R}$ belong to $L^{2}\left( 
\mathbb{R}\right) $, and suppose that its Fourier transform $\widehat{K}$
satisfies the following conditions:%
\begin{gather*}
\widehat{K}\in C^{\infty }\left( \mathbb{R}\setminus \left\{ 0\right\}
\right) ; \\
\sup_{y\in \mathbb{R}\setminus \left\{ 0\right\} }\left( \left| \widehat{K}%
\left( y\right) \right| \max \left\{ 1,\left| y\right| \right\} \right)
<\infty ; \\
\sup_{y\in \mathbb{R}\setminus \left\{ 0\right\} }\left( \left| \frac{d^{n}\,%
\widehat{K}\left( y\right) }{dy^{n}}\right| \max \left\{ \left| y\right|
^{n-1},\left| y\right| ^{n+1}\right\} \right) <\infty \text{, for each }n\in 
\mathbb{N}\text{.}
\end{gather*}%
Suppose that $m\in \mathbb{N}$, and let $d_{m}=2^{1/m}$. Then there is a
positive constant $\gamma _{K,m}$, depending only on $K$ and $m$, such that
for every pair of compactly supported functions $f$ and $g$ belonging to $%
L^{\infty }\left( \mathbb{R}\right) $, for each positive integer $J\geq 2$,
and for each sequence $\left\{ u_{j}\right\} _{j=1}^{J}\subseteq \mathbb{Z}$
such that $u_{1}<u_{2}<\cdots <u_{J}$, we have, in terms of the notation (%
\ref{e1.7}) for dilations,%
\begin{align}
& \left\| \left\{ \sum_{j=1}^{J-1}\sup_{\substack{ k\in \mathbb{Z},  \\ %
u_{j}\leq k<u_{j+1}}}\left| \int_{\mathbb{R}}f\left( x+y\right) g\left(
x-y\right) \left( \left( \delta _{d_{m}^{k}}K\right) \left( y\right) -\left(
\delta _{d_{m}^{u_{j+1}}}K\right) \left( y\right) \right) dy\right|
^{2}\right\} ^{%
{\frac12}%
}\right\| _{L_{x}^{1,\infty }}\hspace{1in}  \label{e2a.0} \\
& \leq \gamma _{K,m}\,J^{1/4}\,\left\| f\right\| _{L^{2}\left( \mathbb{R}%
\right) }\,\left\| g\right\| _{L^{2}\left( \mathbb{R}\right) }\text{.} 
\notag
\end{align}
\end{theorem}

The discretization of Theorem \ref{Dem_Thm_1_4} will be accomplished by
making suitable use of a general discretization tool for maximal functions
(Lemma 3 of \cite{BBCG}, whose statement is reproduced below in Lemma \ref%
{BBCG_Lemma3}). Before embarking on this\ course, we introduce some
auxiliary notions and notation in order to avoid digressions later on. The
closed interval $\left[ -\dfrac{1}{4},\dfrac{1}{4}\right] $ in $\mathbb{R}$
will be designated by $\mathcal{I}$, and for each $n\in \mathbb{Z}$, we
denote by $\mathcal{I}_{n}$ the interval $\mathcal{I}+n=\left[ n-\dfrac{1}{4}%
,n+\dfrac{1}{4}\right] $. For each $\phi \in L^{1}\left( \mathbb{R}\right) $
such that the support of $\phi $ is a subset of $\mathcal{I}$, we define the
linear mapping $P_{\phi }:\mathbb{C}^{\mathbb{Z}}\rightarrow \mathbb{C}^{%
\mathbb{R}}$ by putting%
\begin{eqnarray}
\left( P_{\phi }\left( \left\{ a_{n}\right\} _{n=-\infty }^{\infty }\right)
\right) \left( x\right) &=&\sum_{n\in \mathbb{Z}}a_{n}\,\phi \left(
x-n\right) \text{,}  \label{e2a.1} \\
\text{for all }\left\{ a_{n}\right\} _{n=-\infty }^{\infty } &\in &\mathbb{C}%
^{\mathbb{Z}}\text{, and all }x\in \mathbb{R}\text{.}  \notag
\end{eqnarray}%
When $\phi $ is specialized to be $\chi _{\mathcal{I}}$, the characteristic
function, defined on $\mathbb{R}$, of $\mathcal{I}$, we shall write $P$
rather than $P_{\phi }$. Clearly, if $0<p\leq \infty $, and if $\phi \in
L^{1}\left( \mathbb{R}\right) \bigcap L^{p}\left( \mathbb{R}\right) $ with
support contained in $\mathcal{I}$, then%
\begin{gather}
\left\| P_{\phi }\left( \left\{ a_{n}\right\} _{n=-\infty }^{\infty }\right)
\right\| _{L^{p}\left( \mathbb{R}\right) }=\left\| \phi \right\|
_{L^{p}\left( \mathbb{R}\right) }\left\| \left\{ a_{n}\right\} _{n=-\infty
}^{\infty }\right\| _{\ell ^{p}\left( \mathbb{Z}\right) }\text{,}
\label{e2a.2} \\
\text{for all }\left\{ a_{n}\right\} _{n=-\infty }^{\infty }\in \ell
^{p}\left( \mathbb{Z}\right) \text{.}  \notag
\end{gather}%
Notice also that if $\phi \in L^{1}\left( \mathbb{R}\right) $ is a
non-negative function with support contained in $\mathcal{I}$, if $N\in 
\mathbb{N}$, if, for $1\leq j\leq N$, $a^{\left( j\right) }\equiv \left\{
a_{n}^{\left( j\right) }\right\} _{n=-\infty }^{\infty }$ is a sequence of
real numbers, and if we put%
\begin{equation*}
a_{n}^{\#}=\sup_{1\leq j\leq N}a_{n}^{\left( j\right) }\text{, for all }n\in 
\mathbb{Z}\text{,}
\end{equation*}%
then pointwise on $\mathbb{R}$ we have%
\begin{equation}
P_{\phi }\left( \left\{ a_{n}^{\#}\right\} _{n=-\infty }^{\infty }\right)
=\sup_{1\leq j\leq N}P_{\phi }\left( \left\{ a_{n}^{\left( j\right)
}\right\} _{n=-\infty }^{\infty }\right) \text{.}  \label{e2a.3}
\end{equation}%
For a given $F\in L^{1}\left( \mathbb{R}\right) $, we shall denote by $S_{F}$
the bilinear mapping of $L^{2}\left( \mathbb{R}\right) \times L^{2}\left( 
\mathbb{R}\right) $ into $L^{1}\left( \mathbb{R}\right) $ specified by%
\begin{equation}
\left( S_{F}\left( f,g\right) \right) \left( x\right) =\int_{\mathbb{R}%
}\,f\left( x+y\right) g\left( x-y\right) F\left( y\right) \,dy.
\label{e2a.4}
\end{equation}%
Similarly, for a given $\mathcal{W\equiv }\left\{ \mathcal{W}_{n}\right\}
_{n=-\infty }^{\infty }\in \ell ^{1}\left( \mathbb{Z}\right) $, we define
the bilinear mapping $\mathfrak{S}_{\mathcal{W}}:\ell ^{2}\left( \mathbb{Z}%
\right) \,\times \,\ell ^{2}\left( \mathbb{Z}\right) \rightarrow \ell
^{1}\left( \mathbb{Z}\right) $ by writing for all $\mathfrak{a}\in \,\ell
^{2}\left( \mathbb{Z}\right) $, all $\mathfrak{b}\in \,\ell ^{2}\left( 
\mathbb{Z}\right) $, and all $n\in \mathbb{Z}$, 
\begin{equation}
\left( \mathfrak{S}_{\mathcal{W}}\left( \mathfrak{a},\mathfrak{b}\right)
\right) \left( n\right) =\sum_{k=-\infty }^{\infty }\mathfrak{a}_{n+k}%
\mathfrak{b}_{n-k}\mathcal{W}_{k}\text{.}  \label{e2a.4a}
\end{equation}%
Using the foregoing notation, we can reproduce the statement of Lemma 3 from %
\cite{BBCG} as follows.

\begin{lemma}
\label{BBCG_Lemma3}Let $N\in \mathbb{N}$, and suppose that $\left\{
F_{j}\right\} _{j=1}^{N}\subseteq L^{1}\left( \mathbb{R}\right) $ is such
that for $1\leq j\leq N$, and each $n\in \mathbb{Z}$, the restriction $%
\left( F_{j}\left| \mathcal{I}_{n}\right. \right) $ belongs to $C^{1}\left( 
\mathcal{I}_{n}\right) $, and denote its derivative by $\left( F_{j}\left| 
\mathcal{I}_{n}\right. \right) ^{\prime }$. For each $n\in \mathbb{Z}$, let%
\begin{equation*}
\Lambda _{n}=\sup \left\{ \left| \left( F_{j}\left| \mathcal{I}_{n}\right.
\right) ^{\prime }\left( x\right) \right| :1\leq j\leq N\text{, }x\in 
\mathcal{I}_{n}\right\} \text{,}
\end{equation*}%
and assume that the sequence $\Lambda \equiv \left\{ \Lambda _{n}\right\}
_{n=-\infty }^{\infty }\in \ell ^{1}\left( \mathbb{Z}\right) $. Let%
\begin{gather}
\left( S^{\left( N\right) }\left( f,g\right) \right) \left( x\right) =\sup 
_{\substack{ j\in \mathbb{N},  \\ 1\leq j\leq N}}\left| \left(
S_{F_{j}}\left( f,g\right) \right) \left( x\right) \right| \text{,}
\label{e2a.5} \\
\hspace{0.5in}\text{for all }f,g\in L^{2}\left( \mathbb{R}\right) \text{,
and all }x\in \mathbb{R}\text{,}  \notag
\end{gather}%
and put%
\begin{gather}
\left( \mathfrak{S}^{\left( N\right) }\left( \mathfrak{a},\mathfrak{b}%
\right) \right) \left( m\right) =\sup_{1\leq j\leq N}\left| \sum_{n=-\infty
}^{\infty }\mathfrak{a}_{m+n}\mathfrak{b}_{m-n}F_{j}\left( n\right) \right| 
\text{,}  \label{e2a.6} \\
\text{\ for all }\mathfrak{a},\mathfrak{b}\in \ell ^{2}\left( \mathbb{Z}%
\right) \text{, and all }m\in \mathbb{Z}\text{.}  \notag
\end{gather}%
Further, let $\phi _{0}\geq 0$ and $\phi _{1}\geq 0$ be the functions
defined on $\mathbb{R}$ by writing for each $u\in \mathbb{R}$,%
\begin{eqnarray}
\phi _{0}\left( u\right) &=&2\left( \frac{1}{4}-\left| u\right| \right) \chi
_{\mathcal{I}}\left( u\right) \text{;}  \label{e2a.7} \\
\phi _{1}\left( u\right) &=&\left( \frac{1}{4}-\left| u\right| \right)
^{2}\chi _{\mathcal{I}}\left( u\right) \text{.}  \label{e2a.8}
\end{eqnarray}%
Then for every pair $a,b$ of finitely supported complex-valued sequences
defined on $\mathbb{Z}$, the following inequality holds pointwise on $%
\mathbb{R}$.%
\begin{equation}
P_{\phi _{0}}\left( \mathfrak{S}^{\left( N\right) }\left( a,b\right) \right)
\leq S^{\left( N\right) }\left( Pa,Pb\right) +P_{\phi _{1}}\left( \mathfrak{S%
}_{\Lambda }\left( \left| a\right| ,\left| b\right| \right) \right) \text{.}
\label{e2a.9}
\end{equation}
\end{lemma}

This discussion has paved the way for our discretization of Theorem \ref%
{Dem_Thm_1_4} in the following form. (N.B.-In \S \ref{sec3}, we shall need
to modify the proof of this discretization theorem in order to obtain a
variant of its conclusion for a certain kernel lacking compact support. For
the sake of this subsequent modification we shall keep careful track of the
constants involved in the present proof--particularly as regards the
parameters each of them depends on.)

\begin{theorem}
\label{discrete_Thm1_4}Suppose that $K:$ $\mathbb{R\rightarrow R}$ belongs
to $C^{\infty }\left( \mathbb{R}\right) $ and has compact support. We put%
\begin{eqnarray}
\alpha _{K} &=&\min \left\{ n\in \mathbb{N}:K\left( x\right) =0\text{
whenever }\left| x\right| >n\right\} \text{;}  \label{e2a.9a} \\
\beta _{K} &=&\sup_{x\in \mathbb{R}}\left| K^{\prime }\left( x\right)
\right| \text{.}  \notag
\end{eqnarray}%
Let $m\in \mathbb{N}$, and put $d_{m}=2^{1/m}$. For each $j\in \mathbb{N}$,
let $K_{j,m}\in $ $C^{\infty }\left( \mathbb{R}\right) $ be the compactly
supported function specified by writing (in accordance with the notation for
dilations defined in (\ref{e1.7})) \ $K_{j,m}=\delta _{{\Large d}_{m}^{j}}K$%
, and let $Q_{j,K,m}:\mathbb{C}^{\mathbb{Z}}$ \thinspace $\times \,\mathbb{C}%
^{\mathbb{Z}}\rightarrow \mathbb{C}^{\mathbb{Z}}$ be the bilinear mapping
defined for all $\mathfrak{v}\in \mathbb{C}^{\mathbb{Z}}$ and all $\mathfrak{%
w}\in \mathbb{C}^{\mathbb{Z}}$ by%
\begin{equation}
\left( Q_{j,K,m}\left( \mathfrak{v},\mathfrak{w}\right) \right) \left(
n\right) =\sum_{k=-\infty }^{\infty }\mathfrak{v}\left( n+k\right) \mathfrak{%
w}\left( n-k\right) K_{j,m}\left( k\right) \text{, for all }n\in \mathbb{Z}%
\text{.}  \label{e2a.10}
\end{equation}%
(Notice, in particular, that if $\mathfrak{v}$ and $\mathfrak{w}$ are
finitely supported, then so is the sequence $Q_{{\Large j,K,m}}\left( 
\mathfrak{v},\mathfrak{w}\right) $.) Then for every pair of finitely
supported sequences $a\in \mathbb{C}^{\mathbb{Z}}$ and $b\in \mathbb{C}^{%
\mathbb{Z}}$, for every integer $R\geq 2$, and for each sequence of positive
integers $u_{1}<u_{2}<\cdots <u_{R}$, we have:%
\begin{eqnarray}
&&\left\| \left\{ \sum_{r=1}^{R-1}\sup_{\substack{ j\in \mathbb{N},  \\ %
u_{r}\leq j<u_{r+1}}}\left| Q_{{\Large j,K,m}}\left( a,b\right) -Q_{{\Large u%
}_{r+1}{\Large ,K,m}}\left( a,b\right) \right| ^{2}\right\} ^{%
{\frac12}%
}\right\| _{\ell ^{1,\infty }\left( \mathbb{Z}\right) }  \label{e2a.11} \\
&\leq &\mathfrak{c}_{m}\left( \gamma _{K,m}+\alpha _{K}\beta _{K}\right)
R^{1/4}\,\left\| a\right\| _{\ell ^{2}\left( \mathbb{Z}\right) }\left\|
b\right\| _{\ell ^{2}\left( \mathbb{Z}\right) }\text{,}  \notag
\end{eqnarray}%
where $\mathfrak{c}_{m}$ is a positive constant depending only on $m$, and $%
\gamma _{K,m}$ is the positive constant depending only on $K$ and $m$ that
occurs in (\ref{e2a.0}).
\end{theorem}

\begin{proof}
Until further notice we now fix $r\in \mathbb{N}$ such that $1\leq r\leq R-1$%
. For each $j\in \mathbb{N}$ such that $u_{r}\leq j<u_{r+1}$, define the
compactly supported $C^{\infty }\left( \mathbb{R}\right) $ function $F_{j,r}$
by writing 
\begin{equation}
F_{j,r}=K_{j,m}-K_{u_{r+1},m}\text{.}  \label{e2a.12}
\end{equation}%
Continuing with the notation used in Lemma \ref{BBCG_Lemma3}, we now define
the bilateral sequence $\Lambda ^{\left( r,K,m\right) }\equiv \left\{
\Lambda _{n}^{\left( r,K,m\right) }\right\} _{n=-\infty }^{\infty }$ by
writing for each $n\in \mathbb{Z}$, 
\begin{equation}
\Lambda _{n}^{\left( r,K,m\right) }=\sup \left\{ \left| \left( F_{j,r}\left| 
\mathcal{I}_{n}\right. \right) ^{\prime }\left( x\right) \right| :u_{r}\leq
j<u_{r+1}\text{, }x\in \mathcal{I}_{n}\right\} \text{.}  \label{e2a.13}
\end{equation}%
Notice that $\left\{ \Lambda _{n}^{\left( r,K,m\right) }\right\} _{n=-\infty
}^{\infty }$ is finitely supported, since $K$ has compact support, and so we
can apply the conclusion (\ref{e2a.9}) of Lemma \ref{BBCG_Lemma3} to the
finite sequence of functions $\left\{ F_{j,r}\right\} _{j=u_{r}}^{u_{r+1}-1}$
and its corresponding sequence $\left\{ \Lambda _{n}^{\left( r,K,m\right)
}\right\} _{n=-\infty }^{\infty }$ specified by (\ref{e2a.13}). This
furnishes our present circumstances with the following pointwise inequality
on $\mathbb{R}$:%
\begin{eqnarray}
&&P_{\phi _{0}}\left( \sup_{\substack{ j\in \mathbb{N},  \\ u_{r}\leq
j<u_{r+1}}}\left| Q_{{\Large j,K,m}}\left( a,b\right) -Q_{{\Large u}_{r+1}%
{\Large ,K,m}}\left( a,b\right) \right| \right)  \label{e2a.14} \\
&\leq &\sup_{\substack{ j\in \mathbb{N}  \\ u_{r}\leq j<u_{r+1}}}\left|
S_{F_{j,r}}\left( Pa,Pb\right) \right| +P_{\phi _{1}}\left( \mathfrak{S}_{%
{\Large \Lambda }^{\left( r,K,m\right) }}\left( \left| a\right| ,\left|
b\right| \right) \right) \text{,}  \notag
\end{eqnarray}%
where the bilinear form $\mathfrak{S}_{{\Large \Lambda }^{\left(
r,K,m\right) }}:\ell ^{2}\left( \mathbb{Z}\right) \,\times \,\ell ^{2}\left( 
\mathbb{Z}\right) \rightarrow \ell ^{1}\left( \mathbb{Z}\right) $ is defined
in accord with (\ref{e2a.4a}), and thus satisfies%
\begin{equation}
\left( \mathfrak{S}_{{\Large \Lambda }^{\left( r,K,m\right) }}\left( \left|
a\right| ,\left| b\right| \right) \right) =\left\{ \sum_{k=-\infty }^{\infty
}\left| a_{n+k}b_{n-k}\right| \Lambda _{k}^{\left( r,K,m\right) }\right\}
_{n=-\infty }^{\infty }\text{.}  \label{e2a.14a}
\end{equation}%
Before letting $r\in \mathbb{N}$ run through all values $1\leq r\leq R-1$ as
on the left side of (\ref{e2a.11}), we first estimate $\left\| \Lambda
^{\left( r,K,m\right) }\right\| _{\ell ^{1}\left( \mathbb{Z}\right) }$ for a
fixed $r$ in this range. For $j\in \mathbb{N}$ satisfying $u_{r}\leq
j<u_{r+1}$, we have $F_{j,r}\left( x\right) =0$, whenever $\left| x\right| $ 
$>\alpha _{K}d_{m}^{u_{r+1}}$, and consequently%
\begin{equation}
\Lambda _{n}^{\left( r,K,m\right) }=0\text{, for all }n\in \mathbb{Z}\text{
such that }\left| n\right| >\alpha _{K}d_{m}^{u_{r+1}}+\frac{1}{4}\text{.}
\label{e2a.17}
\end{equation}%
Suppose that $s\in \mathbb{N}$ with $u_{r}\leq s<u_{r+1}$. Then for all $%
j\in \mathbb{N}$ such that $s\leq j<u_{r+1}$, we have for all $x\in \mathbb{R%
}$, 
\begin{eqnarray}
&&\left| K_{j,m}^{\prime }\left( x\right) -K_{u_{r+1},m}^{\prime }\left(
x\right) \right|  \label{e2a.19} \\
&=&\left| \frac{1}{d_{m}^{2j}}K^{\prime }\left( \frac{x}{d_{m}^{j}}\right) -%
\frac{1}{d_{m}^{2u_{r+1}}}K^{\prime }\left( \frac{x}{d_{m}^{u_{r+1}}}\right)
\right|  \notag \\
&\leq &\beta _{K}\left( \frac{1}{d_{m}^{2s}}+\frac{1}{d_{m}^{2u_{r+1}}}%
\right) \text{.}  \notag
\end{eqnarray}%
Moreover, if $n\in \mathbb{Z}$ and satisfies%
\begin{equation}
\left| n\right| \geq \alpha _{K}d_{m}^{s}+\frac{1}{4}\text{,}  \label{e2a.18}
\end{equation}%
then for all $x\in \mathcal{I}_{n}$ and for all $j\in \mathbb{N}$ such that $%
j<s$,%
\begin{equation}
K_{j,m}^{\prime }\left( x\right) =0\text{.}  \label{e2a.20a}
\end{equation}%
Combining this with (\ref{e2a.19}) we infer that for each $s\in \mathbb{N}$
satisfying $u_{r}\leq s<u_{r+1}$, the following estimate is valid: 
\begin{equation}
\Lambda _{n}^{\left( r,K,m\right) }\leq \beta _{K}\left( \frac{1}{d_{m}^{2s}}%
+\frac{1}{d_{m}^{2u_{r+1}}}\right) \text{, for all }n\in \mathbb{Z}\text{
such that }\left| n\right| \geq \alpha _{K}d_{m}^{s}+\frac{1}{4}\text{.}
\label{e2a.20}
\end{equation}%
Consequently, we see with the aid of (\ref{e2a.17}) that%
\begin{eqnarray}
&&\sum \left\{ \Lambda _{n}^{\left( r,K,m\right) }:n\in \mathbb{Z}\text{, }%
\left| n\right| \geq \alpha _{K}d_{m}^{u_{r}}+\frac{1}{4}\right\}
\label{e2a.21} \\
&\leq &\sum_{\substack{ s\in \mathbb{N}  \\ u_{r}\leq s<u_{r+1}}}\sum
\left\{ \Lambda _{n}^{\left( r,K,m\right) }:n\in \mathbb{Z}\text{, }\alpha
_{K}d_{m}^{s+1}+\frac{1}{4}\geq \left| n\right| \geq \alpha _{K}d_{m}^{s}+%
\frac{1}{4}\right\}  \notag \\
&\leq &2\alpha _{K}\beta _{K}d_{m}\sum_{_{_{\substack{ s\in \mathbb{N}  \\ %
u_{r}\leq s<u_{r+1}}}}}\left( \frac{1}{d_{m}^{2s}}+\frac{1}{d_{m}^{2u_{r+1}}}%
\right) d_{m}^{s}  \notag \\
&\leq &\left( \frac{2\alpha _{K}\beta _{K}d_{m}^{2}}{d_{m}-1}\right) \frac{1%
}{d_{m}^{u_{r}}}+2\alpha _{K}\beta _{K}d_{m}\frac{u_{r+1}}{d_{m}^{u_{r+1}}}%
\text{.}  \notag
\end{eqnarray}%
By specializing the value of $s$ in (\ref{e2a.19}) to be $u_{r}$, we see
that for all $n\in \mathbb{Z}$,%
\begin{equation*}
\Lambda _{n}^{\left( r,K,m\right) }\leq \beta _{K}\left( \frac{1}{%
d_{m}^{2u_{r}}}+\frac{1}{d_{m}^{2u_{r+1}}}\right) \text{.}
\end{equation*}%
Hence%
\begin{eqnarray}
&&\sum \left\{ \Lambda _{n}^{\left( r,K,m\right) }:n\in \mathbb{Z}\text{, }%
\left| n\right| \leq \alpha _{K}d_{m}^{u_{r}}+\frac{1}{4}\right\}
\label{e2a.22} \\
&\leq &\beta _{K}\left( \frac{1}{d_{m}^{2u_{r}}}+\frac{1}{d_{m}^{2u_{r+1}}}%
\right) \left( 2\alpha _{K}d_{m}^{u_{r}}+\frac{3}{2}\right)  \notag \\
&\leq &2\alpha _{K}\beta _{K}\left( \frac{1}{d_{m}^{u_{r}}}+\frac{1}{%
d_{m}^{u_{r+1}}}\right) +\frac{3}{2}\beta _{K}\left( \frac{1}{d_{m}^{2u_{r}}}%
+\frac{1}{d_{m}^{2u_{r+1}}}\right) \text{.}  \notag
\end{eqnarray}%
Upon combining (\ref{e2a.21}) and (\ref{e2a.22}), we deduce that for each $%
r\in \mathbb{N}$ such that $1\leq r\leq R-1$, 
\begin{equation}
\left\| \Lambda ^{\left( r,K,m\right) }\right\| _{\ell ^{1}\left( \mathbb{Z}%
\right) }\leq \alpha _{K}\beta _{K}C_{m}\left( \frac{1}{d_{m}^{u_{r}}}+\frac{%
u_{r+1}}{d_{m}^{u_{r+1}}}\right) \text{.}  \label{e2a.23}
\end{equation}

We next square both sides of (\ref{e2a.14}). In view of the definitions of $%
\phi _{0}$ and $\phi _{1}$ in (\ref{e2a.7}) and (\ref{e2a.8}) this gives us
the following inequality, valid pointwise on $\mathbb{R}$.%
\begin{eqnarray}
&&4P_{\phi _{1}}\left( \sup_{\substack{ j\in \mathbb{N},  \\ u_{r}\leq
j<u_{r+1}}}\left| Q_{{\Large j,K,m}}\left( a,b\right) -Q_{{\Large u}_{r+1}%
{\Large ,K,m}}\left( a,b\right) \right| ^{2}\right)  \label{e2a.24} \\
&\leq &\left\{ \sup_{\substack{ j\in \mathbb{N}  \\ u_{r}\leq j<u_{r+1}}}%
\left| S_{F_{j,r}}\left( Pa,Pb\right) \right| +P_{\phi _{1}}\left( \mathfrak{%
S}_{{\Large \Lambda }^{\left( r,K,m\right) }}\left( \left| a\right| ,\left|
b\right| \right) \right) \right\} ^{2}\text{.}  \notag
\end{eqnarray}%
After summing this inequality for $1\leq r\leq R-1$, we deduce with the aid
of Minkowski's inequality for $\ell ^{2}$ that the following holds pointwise
on $\mathbb{R}$.%
\begin{eqnarray}
&&P_{\phi _{0}}\left( \left\{ \sum_{r=1}^{R-1}\left( \sup_{\substack{ j\in 
\mathbb{N},  \\ u_{r}\leq j<u_{r+1}}}\left| Q_{{\Large j,K,m}}\left(
a,b\right) -Q_{{\Large u}_{r+1}{\Large ,K,m}}\left( a,b\right) \right|
^{2}\right) \right\} ^{%
{\frac12}%
}\right)  \label{e2a.25} \\
&\leq &\left\{ \sum_{r=1}^{R-1}\sup_{\substack{ j\in \mathbb{N}  \\ %
u_{r}\leq j<u_{r+1}}}\left| S_{F_{j,r}}\left( Pa,Pb\right) \right|
^{2}\right\} ^{%
{\frac12}%
}+P_{\phi _{1}}\left( \sum_{r=1}^{R-1}\mathfrak{S}_{{\Large \Lambda }%
^{\left( r,K,m\right) }}\left( \left| a\right| ,\left| b\right| \right)
\right) \text{.}  \notag
\end{eqnarray}%
For each $n\in \mathbb{Z}$, $\phi _{0}\left( x-n\right) \geq 1/4$ provided $%
\left| x-n\right| \leq 1/8$. It follows readily that, in terms of
distribution functions (taken with respect to Lebesgue measure $\mathfrak{m}%
_{\mathbb{R}}$ on $\mathbb{R}$ and counting measure $\mathfrak{m}_{\mathbb{Z}%
}$ on $\mathbb{Z}$), we have for each positive real number $y$,%
\begin{gather}
\frac{1}{4}\lambda \left( \left\{ \sum_{r=1}^{R-1}\left( \sup_{\substack{ %
j\in \mathbb{N},  \\ u_{r}\leq j<u_{r+1}}}\left| Q_{{\Large j,K,m}}\left(
a,b\right) -Q_{{\Large u}_{r+1}{\Large ,K,m}}\left( a,b\right) \right|
^{2}\right) \right\} ^{%
{\frac12}%
},\mathfrak{m}_{\mathbb{Z}};y\right)  \label{e2a.26} \\
\leq \lambda \left( P_{\phi _{0}}\left( \left\{ \sum_{r=1}^{R-1}\left( \sup 
_{\substack{ j\in \mathbb{N},  \\ u_{r}\leq j<u_{r+1}}}\left| Q_{{\Large %
j,K,m}}\left( a,b\right) -Q_{{\Large u}_{r+1}{\Large ,K,m}}\left( a,b\right)
\right| ^{2}\right) \right\} ^{%
{\frac12}%
}\right) ,\mathfrak{m}_{\mathbb{R}};\frac{y}{4}\right) \text{.}  \notag
\end{gather}%
From (\ref{e2a.25}) and (\ref{e2a.26}) we see at once that%
\begin{gather}
\frac{1}{4}\lambda \left( \left\{ \sum_{r=1}^{R-1}\left( \sup_{\substack{ %
j\in \mathbb{N},  \\ u_{r}\leq j<u_{r+1}}}\left| Q_{{\Large j,K,m}}\left(
a,b\right) -Q_{{\Large u}_{r+1}{\Large ,K,m}}\left( a,b\right) \right|
^{2}\right) \right\} ^{%
{\frac12}%
},\mathfrak{m}_{\mathbb{Z}};y\right)  \label{e2a.28} \\
\leq \lambda \left( \left\{ \sum_{r=1}^{R-1}\sup_{\substack{ j\in \mathbb{N} 
\\ u_{r}\leq j<u_{r+1}}}\left| S_{F_{j,r}}\left( Pa,Pb\right) \right|
^{2}\right\} ^{%
{\frac12}%
},\mathfrak{m}_{\mathbb{R}};\frac{y}{8}\right)  \notag \\
\hspace{0.3in}+\lambda \left( P_{\phi _{1}}\left( \sum_{r=1}^{R-1}\mathfrak{S%
}_{{\Large \Lambda }^{\left( r,K,m\right) }}\left( \left| a\right| ,\left|
b\right| \right) \right) ,\mathfrak{m}_{\mathbb{R}};\frac{y}{8}\right) \text{%
.}  \notag
\end{gather}%
Moreover, an application of Theorem \ref{Dem_Thm_1_4} shows that%
\begin{eqnarray}
&&\lambda \left( \left\{ \sum_{r=1}^{R-1}\sup_{\substack{ j\in \mathbb{N} 
\\ u_{r}\leq j<u_{r+1}}}\left| S_{F_{j},r}\left( Pa,Pb\right) \right|
^{2}\right\} ^{%
{\frac12}%
},\mathfrak{m}_{\mathbb{R}};\frac{y}{8}\right)  \label{e2a.29} \\
&\leq &\frac{4\gamma _{K,m}\,R^{1/4}\left\| a\right\| _{\ell ^{2}\left( 
\mathbb{Z}\right) }\,\left\| b\right\| _{\ell ^{2}\left( \mathbb{Z}\right) }%
}{y}\text{.}  \notag
\end{eqnarray}
Applying Chebychev's inequality to the function $P_{\phi _{1}}\left(
\sum_{r=1}^{R-1}\mathfrak{S}_{{\Large \Lambda }^{\left( r,K,m\right)
}}\left( \left| a\right| ,\left| b\right| \right) \right) $, we find with
the aid of (\ref{e2a.23}) that 
\begin{equation}
\lambda \left( P_{\phi _{1}}\left( \sum_{r=1}^{R-1}\mathfrak{S}_{{\Large %
\Lambda }^{\left( r,K,m\right) }}\left( \left| a\right| ,\left| b\right|
\right) \right) ,\mathfrak{m}_{\mathbb{R}};\frac{y}{8}\right) \leq \frac{%
\alpha _{K}\beta _{K}C_{m}\left\| a\right\| _{\ell ^{2}\left( \mathbb{Z}%
\right) }\,\left\| b\right\| _{\ell ^{2}\left( \mathbb{Z}\right) }}{y}\text{.%
}  \label{e2a.30}
\end{equation}%
The desired conclusion (\ref{e2a.11}) is an immediate consequence of (\ref%
{e2a.28}), (\ref{e2a.29}), and (\ref{e2a.30}).
\end{proof}

The next theorem and its included corollary of Theorem \ref{discrete_Thm1_4}
(see Corollary \ref{Transferred_Oscillation} below) furnish key applications
of Lemma \ref{Demeter_Lemma3_1} by using the ``isometric'' transformation $U$
in the hypotheses of Theorem \ref{Main_Result} to transfer\ discrete
oscillation estimates such as (\ref{e2a.11}) of Theorem \ref{discrete_Thm1_4}
to the setting of the arbitrary sigma-finite measure space $\left( \Omega
,\mu \right) $. (Compare the transference reasoning in \cite{Dem} and \cite%
{DTT} that targeted averages defined by the invertible measure-preserving
point transformations of finite measure spaces.) The notation of Proposition %
\ref{structure_of_U} will now be in effect, and it will be convenient to
observe that since each $\Phi _{k}$, $k\in \mathbb{Z}$, is multiplicative on 
$\mathcal{A}\left( \mu \right) $, it follows from Proposition \ref%
{structure_of_U}-(jj), together with (\ref{e2.4}) and (\ref{e2.5}), that for
all $F\in \mathcal{A}\left( \mu \right) $, $G\in \mathcal{A}\left( \mu
\right) $, $n_{1}\in \mathbb{Z}$, and $n_{2}\in \mathbb{Z}$, we have $\mu $%
-a.e. on $\Omega $,%
\begin{eqnarray}
\Phi _{k}\left( \left( U^{n_{1}}F\right) \left( U^{n_{2}}G\right) \right)
&=&\Phi _{k}\left( h_{n_{1}}\right) \,\Phi _{k+n_{1}}\left( F\right) \,\Phi
_{k}\left( h_{n_{2}}\right) \,\Phi _{k+n_{2}}\left( G\right)  \label{e2a.31}
\\
&=&h_{k}^{-2}\,\left( U^{k+n_{1}}F\right) \,\left( U^{k+n_{2}}G\right) \text{%
.}  \notag
\end{eqnarray}

\begin{theorem}
\label{gen_discrete_tran}Suppose that $\left( \Omega ,\mu \right) $ is a
sigma-finite measure space, and let $U$ be a bijective linear mapping of $%
\mathcal{A}\left( \mu \right) $ onto $\mathcal{A}\left( \mu \right) $ such
that conditions (i) and (ii) in the hypotheses of Theorem \ref{Main_Result}
hold. For each $j\in \mathbb{N}$, let $\mathfrak{s}_{j}:\mathbb{Z}%
\rightarrow \mathbb{C}$ be finitely supported, and define the bilinear
mappings $T_{j}:\mathbb{C}^{\mathbb{Z}}\,\times \,\mathbb{C}^{\mathbb{Z}%
}\rightarrow \mathbb{C}^{\mathbb{Z}}$ and $\mathfrak{T}_{j}:\mathcal{A}%
\left( \mu \right) \,\times \,\mathcal{A}\left( \mu \right) \rightarrow 
\mathcal{A}\left( \mu \right) $ as follows. 
\begin{gather*}
\left( T_{j}\left( \mathfrak{v,w}\right) \right) \left( k\right)
=\sum_{n=-\infty }^{\infty }\mathfrak{v}\left( k+n\right) \,\mathfrak{w}%
\left( k-n\right) \,\mathfrak{s}_{j}\left( n\right) \text{,} \\
\hspace{0.3in}\hspace{0.3in}\hspace{0.3in}\hspace{0.15in}\text{for each }%
\mathfrak{v}\in \mathbb{C}^{\mathbb{Z}}\text{, each }\mathfrak{w}\in \mathbb{%
C}^{\mathbb{Z}}\text{, and all }k\in \mathbb{Z}\text{.} \\
\mathfrak{T}_{j}\left( F,G\right) =\sum_{n=-\infty }^{\infty }\left(
U^{n}F\right) \,\left( U^{-n}G\right) \mathfrak{s}_{j}\left( n\right) \text{,%
} \\
\hspace{0.6in}\text{for all }F\in \mathcal{A}\left( \mu \right) \text{, and
all }G\in \mathcal{A}\left( \mu \right) \text{.}
\end{gather*}%
Suppose that there is a positive real constant $\zeta $ such that for every
pair of finitely supported sequences $a\in \mathbb{C}^{\mathbb{Z}}$ and $%
b\in \mathbb{C}^{\mathbb{Z}}$, for every integer $R\geq 2$, and for each
sequence of positive integers $u_{1}<u_{2}<\cdots <u_{R}$, we have:%
\begin{eqnarray}
&&\left\| \left\{ \sum_{r=1}^{R-1}\sup_{\substack{ j\in \mathbb{N},  \\ %
u_{r}\leq j<u_{r+1}}}\left| T_{j}\left( a,b\right) -T_{u_{r+1}}\left(
a,b\right) \right| ^{2}\right\} ^{%
{\frac12}%
}\right\| _{\ell ^{1,\infty }\left( \mathbb{Z}\right) }  \label{e2a.b1} \\
&\leq &\zeta R^{1/4}\,\left\| a\right\| _{\ell ^{2}\left( \mathbb{Z}\right)
}\left\| b\right\| _{\ell ^{2}\left( \mathbb{Z}\right) }\text{.}  \notag
\end{eqnarray}%
Then for every $f\in L^{2}\left( \mu \right) $, every $g\in L^{2}\left( \mu
\right) $, every integer $R\geq 2$, and each sequence of positive integers $%
u_{1}<u_{2}<\cdots <u_{R}$, we have:%
\begin{eqnarray}
&&\left\| \left\{ \sum_{r=1}^{R-1}\sup_{\substack{ j\in \mathbb{N},  \\ %
u_{r}\leq j<u_{r+1}}}\left| \mathfrak{T}_{j}\left( f,g\right) -\mathfrak{T}%
_{u_{r+1}}\left( f,g\right) \right| ^{2}\right\} ^{1/2}\right\|
_{L^{1,\infty }\left( \mu \right) }  \label{e2a.b2} \\
&\leq &\zeta R^{1/4}\,\left\| f\right\| _{L^{2}\left( \mu \right) }\left\|
g\right\| _{L^{2}\left( \mu \right) }\text{.}  \notag
\end{eqnarray}%
Hence by Lemma \ref{Demeter_Lemma3_1}, for every $f\in L^{2}\left( \mu
\right) $ and every $g\in L^{2}\left( \mu \right) $ the sequence $\left\{ 
\mathfrak{T}_{j}\left( f,g\right) \right\} _{j=1}^{\infty }$ converges $\mu $%
-a.e. on $\Omega $ to a corresponding function belonging to $\mathcal{A}%
\left( \mu \right) $.
\end{theorem}

\begin{proof}
For convenience, put%
\begin{equation}
\Delta =\left\{ \sum_{r=1}^{R-1}\sup_{\substack{ j\in \mathbb{N},  \\ %
u_{r}\leq j<u_{r+1}}}\left| \mathfrak{T}_{j}\left( f,g\right) -\mathfrak{T}%
_{u_{r+1}}\left( f,g\right) \right| ^{2}\right\} ^{1/2}\text{.}
\label{e2a.34}
\end{equation}%
Applying (\ref{e2.6}) to the distribution function of $\Delta $, we see that
for each $k\in \mathbb{Z}$, and each real number $y>0$,%
\begin{equation}
\lambda \left( \Delta ,\mu ;y\right) =\lambda \left( \Phi _{k}\left( \Delta
\right) ,\mu ;y\right) \text{.}  \label{e2a.35}
\end{equation}%
From the properties of $\Phi _{k}$ in (\ref{e2.2}), (\ref{e2.3}), and (\ref%
{e2.5a}) we deduce that%
\begin{eqnarray}
&&\Phi _{k}\left( \Delta \right)  \label{e2a.36} \\
&=&\left\{ \sum_{r=1}^{R-1}\sup_{\substack{ j\in \mathbb{N},  \\ u_{r}\leq
j<u_{r+1}}}\left| \Phi _{k}\left( \mathfrak{T}_{j}\left( f,g\right) \right)
-\Phi _{k}\left( \mathfrak{T}_{u_{r+1}}\left( f,g\right) \right) \right|
^{2}\right\} ^{1/2}  \notag \\
&=&\left\{ \sum_{r=1}^{R-1}\sup_{\substack{ j\in \mathbb{N},  \\ u_{r}\leq
j<u_{r+1}}}\left| h_{k}^{2}\Phi _{k}\left( \mathfrak{T}_{j}\left( f,g\right)
\right) -h_{k}^{2}\Phi _{k}\left( \mathfrak{T}_{u_{r+1}}\left( f,g\right)
\right) \right| ^{2}\right\} ^{1/2}\text{.}  \notag
\end{eqnarray}%
For every $\nu \in \mathbb{N}$, we have by virtue of (\ref{e2a.31}),%
\begin{equation}
h_{k}^{2}\Phi _{k}\left( \mathfrak{T}_{\upsilon }\left( f,g\right) \right)
=\sum_{n=-\infty }^{\infty }\left( U^{k+n}f\right) \,\left( U^{k-n}g\right) 
\mathfrak{s}_{\upsilon }\left( n\right) \text{.}  \label{e2a.37}
\end{equation}%
Using (\ref{e2a.37}) to substitute in (\ref{e2a.36}), we find that for each $%
k\in \mathbb{Z}$, 
\begin{eqnarray}
\Phi _{k}\left( \Delta \right) &=&\left\{ \sum_{r=1}^{R-1}\sup_{\substack{ %
j\in \mathbb{N},  \\ u_{r}\leq j<u_{r+1}}}\left| \sum_{n=-\infty }^{\infty
}\left( U^{k+n}f\right) \,\left( U^{k-n}g\right) \mathfrak{s}_{j}\left(
n\right) -\right. \right.  \label{e2a.38} \\
&&\hspace{0.5in}\left. \left. \sum_{n=-\infty }^{\infty }\left(
U^{k+n}f\right) \,\left( U^{k-n}g\right) \mathfrak{s}_{u_{r+1}}\left(
n\right) \right| ^{2}\right\} ^{1/2}\text{.}  \notag
\end{eqnarray}%
Now let $N_{0}$ be the least positive integer $N$ such that $\mathfrak{s}%
_{j}\left( n\right) =0$ whenever $1\leq j\leq u_{R}$ and $\left| n\right| >N$%
, temporarily fix an arbitrary $L\in \mathbb{N}$, and let $\mathfrak{C}%
_{L,N_{0}}$ denote the characteristic function, defined on $\mathbb{Z}$, of%
\begin{equation*}
\left\{ n\in \mathbb{Z}:\left| n\right| \leq L+N_{0}\right\} \text{.}
\end{equation*}%
For each $x\in \Omega $, we define the finitely supported sequences $\phi
_{x}:\mathbb{Z\rightarrow C}$ and $\psi _{x}:\mathbb{Z\rightarrow C}$ by
writing for each $n\in \mathbb{Z}$,%
\begin{eqnarray*}
\phi _{x}\left( n\right) &=&\mathfrak{C}_{L,N_{0}}\left( n\right) \,\left(
\left( U^{n}f\right) \left( x\right) \right) \text{;} \\
\psi _{x}\left( n\right) &=&\mathfrak{C}_{L,N_{0}}\left( n\right) \,\,\left(
\left( U^{n}g\right) \left( x\right) \right) \text{.}
\end{eqnarray*}%
In terms of this notation, we can use (\ref{e2a.38}) to write for each $k\in 
\mathbb{Z}$ such that $-L\leq k\leq L$, and for each $x\in \Omega $,%
\begin{eqnarray}
&&\left( \Phi _{k}\left( \Delta \right) \right) \left( x\right)
\label{e2a.38a} \\
&=&\left\{ \sum_{r=1}^{R-1}\sup_{\substack{ j\in \mathbb{N},  \\ u_{r}\leq
j<u_{r+1}}}\left| \left( T_{j}\left( \phi _{x},\psi _{x}\right) \right)
\left( k\right) -\left( T_{u_{r+1}}\left( \phi _{x},\psi _{x}\right) \right)
\left( k\right) \right| ^{2}\right\} ^{1/2}\text{.}  \notag
\end{eqnarray}%
From (\ref{e2a.35}) we have for each real number $y>0$,%
\begin{equation}
\left( 2L+1\right) \lambda \left( \Delta ,\mu ;y\right)
=\sum_{k=-L}^{L}\lambda \left( \Phi _{k}\left( \Delta \right) ,\mu ;y\right) 
\text{.}  \label{e2a.39}
\end{equation}%
Temporarily fix an arbitrary positive real number $y$, and for each $k\in 
\mathbb{Z}$ with $-L\leq k\leq L$, denote by $\chi _{k}$ the characteristic
function, defined on $\Omega $, of the set $E_{k}$ specified by%
\begin{equation*}
E_{k}=\left\{ x\in \Omega :\left( \Phi _{k}\left( \Delta \right) \right)
\left( x\right) >y\right\} \text{.}
\end{equation*}%
This permits us to rewrite (\ref{e2a.39}) in the form%
\begin{equation}
\left( 2L+1\right) \lambda \left( \Delta ,\mu ;y\right) =\int_{\Omega
}\,\left( \sum_{k=-L}^{L}\chi _{k}\left( x\right) \right) \,d\mu \left(
x\right) \text{.}  \label{e2a.40}
\end{equation}%
With the aid of (\ref{e2a.38a}) we see that at each $x\in \Omega $, the
integrand in (\ref{e2a.40}) can be expressed in terms of counting measure $%
\mathfrak{m}_{\mathbb{Z}}$ on $\mathbb{Z}$ by:%
\begin{align*}
& \sum_{k=-L}^{L}\chi _{k}\left( x\right) \\
& =\mathfrak{m}_{\mathbb{Z}}\left\{ k\in \mathbb{Z}:-L\leq k\leq L,\text{
and }x\in E_{k}\right\} \\
& \leq \mathfrak{m}_{\mathbb{Z}}\left\{ k\in \mathbb{Z}:\left\{
\sum_{r=1}^{R-1}\sup_{\substack{ j\in \mathbb{N},  \\ u_{r}\leq j<u_{r+1}}}%
\left| \left( T_{j}\left( \phi _{x},\psi _{x}\right) \right) \left( k\right)
-\left( T_{u_{r+1}}\left( \phi _{x},\psi _{x}\right) \right) \left( k\right)
\right| ^{2}\right\} ^{1/2}>y\right\} \text{.}
\end{align*}%
Application to this of the hypothesis (\ref{e2a.b1}) shows that for each $%
x\in \Omega $,%
\begin{eqnarray*}
&&\sum_{k=-L}^{L}\chi _{k}\left( x\right) \\
&\leq &\frac{\zeta R^{1/4}\left\| \phi _{x}\right\| _{\ell ^{2}\left( 
\mathbb{Z}\right) }\,\left\| \psi _{x}\right\| _{\ell ^{2}\left( \mathbb{Z}%
\right) }}{y} \\
&=&\frac{\zeta R^{1/4}}{y}\left\{ \sum_{n=-L-N_{0}}^{L+N_{0}}\left| \left(
U^{n}f\right) \left( x\right) \right| ^{2}\right\} ^{1/2}\,\left\{
\sum_{n=-L-N_{0}}^{L+N_{0}}\left| \left( U^{n}g\right) \left( x\right)
\right| ^{2}\right\} ^{1/2}\text{.}
\end{eqnarray*}%
Using this on the right of (\ref{e2a.40}) and then invoking Cauchy-Schwarz,
we find, since $U\left| L^{2}\left( \mu \right) \right. $ \ is a surjective
linear isometry, that:%
\begin{eqnarray}
&&\lambda \left( \Delta ,\mu ;y\right)  \label{e2a.41} \\
&\leq &\frac{\zeta R^{1/4}}{y\left( 2L+1\right) }\int_{\Omega }\,\left\{
\sum_{n=-L-N_{0}}^{L+N_{0}}\left| \left( U^{n}f\right) \right| ^{2}\right\}
^{1/2}\,\left\{ \sum_{n=-L-N_{0}}^{L+N_{0}}\left| \left( U^{n}g\right)
\right| ^{2}\right\} ^{1/2}\,d\mu  \notag \\
&\leq &\frac{\zeta R^{1/4}}{y}\left( \frac{2L+2N_{0}+1}{2L+1}\right) \left\|
f\right\| _{L^{2}\left( \mu \right) }\,\left\| g\right\| _{L^{2}\left( \mu
\right) }.  \notag
\end{eqnarray}%
In view of the definition of $\Delta $ in (\ref{e2a.34}), we can immediately
arrive at (\ref{e2a.b2}) by letting $L\rightarrow \infty $ on the right of (%
\ref{e2a.41}).
\end{proof}

The following corollary results directly from Theorem \ref{discrete_Thm1_4}
and Theorem \ref{gen_discrete_tran}.

\begin{corollary}
\label{Transferred_Oscillation}Suppose that $\left( \Omega ,\mu \right) $ is
a sigma-finite measure space, and let $U$ be a bijective linear mapping of $%
\mathcal{A}\left( \mu \right) $ onto $\mathcal{A}\left( \mu \right) $ such
that conditions (i) and (ii) in the hypotheses of Theorem \ref{Main_Result}
hold. Suppose that $K:$ $\mathbb{R\rightarrow R}$ belongs to $C^{\infty
}\left( \mathbb{R}\right) $ and has compact support. Let $m\in \mathbb{N}$,
and put $d_{m}=2^{1/m}$. For each $j\in \mathbb{N}$, let $K_{j,m}\in $ $%
C^{\infty }\left( \mathbb{R}\right) $ be the compactly supported function
given by $K_{j,m}=\delta _{{\Large d}_{m}^{j}}K$, and define the bilinear
mapping $\mathfrak{A}_{j,K,m,U}$ of $\mathcal{A}\left( \mu \right) \,\times
\,\mathcal{A}\left( \mu \right) $ into $\mathcal{A}\left( \mu \right) $ by
writing for all $F\in \mathcal{A}\left( \mu \right) $, and all $G\in 
\mathcal{A}\left( \mu \right) $,%
\begin{eqnarray}
\mathfrak{A}_{j,K,m,U}\left( F,G\right) &=&\sum_{n=-\infty }^{\infty }\left(
U^{n}F\right) \,\left( U^{-n}G\right) K_{j,m}\left( n\right)  \label{e2a.32}
\\
&=&\frac{1}{d_{m}^{j}}\,\sum_{n=-\infty }^{\infty }\left( U^{n}F\right)
\,\left( U^{-n}G\right) K\left( \frac{n}{d_{m}^{j}}\right) \text{.}  \notag
\end{eqnarray}%
Then for every $f\in L^{2}\left( \mu \right) $, every $g\in L^{2}\left( \mu
\right) $, every integer $R\geq 2$, and each sequence of positive integers $%
u_{1}<u_{2}<\cdots <u_{R}$, we have:%
\begin{eqnarray}
&&\left\| \left\{ \sum_{r=1}^{R-1}\sup_{\substack{ j\in \mathbb{N},  \\ %
u_{r}\leq j<u_{r+1}}}\left| \mathfrak{A}_{j,K,m,U}\left( f,g\right) -%
\mathfrak{A}_{u_{r+1},K,m,U}\left( f,g\right) \right| ^{2}\right\}
^{1/2}\right\| _{L^{1,\infty }\left( \mu \right) }  \label{e2a.33} \\
&\leq &\Gamma _{K,m}R^{1/4}\,\left\| f\right\| _{L^{2}\left( \mu \right)
}\left\| g\right\| _{L^{2}\left( \mu \right) }\text{,}  \notag
\end{eqnarray}%
where $\Gamma _{K,m}$ denotes the constant $\mathfrak{c}_{m}\left( \gamma
_{K,m}+\alpha _{K}\beta _{K}\right) $ that occurs in (\ref{e2a.11}). Hence
for every $f\in L^{2}\left( \mu \right) $ and every $g\in L^{2}\left( \mu
\right) $, the sequence $\left\{ \mathfrak{A}_{j,K,m,U}\left( f,g\right)
\right\} _{j=1}^{\infty }$ converges $\mu $-a.e. on $\Omega $ to a
corresponding function belonging to $\mathcal{A}\left( \mu \right) $.
\end{corollary}

\section{Proof of Theorem \ref{Main_Result}\label{sec3}}

\noindent In view of Theorem \ref{bisublinear_maxml} and dominated
convergence, we need only demonstrate the conclusions of Theorem \ref%
{Main_Result} regarding the existence of $\mu $-a.e. limits. By Theorem \ref%
{bisublinear_maxml} and the Multilinear Banach Principle (see, e.g.,
Proposition 1 of \cite{BBCG} for the latter), it suffices for the
demonstration of Theorem \ref{Main_Result} to show that each of the
sequences $\left\{ A_{k,U}\left( f,g\right) \right\} _{k=1}^{\infty }$ and $%
\left\{ H_{k,U}\left( f,g\right) \right\} _{k=1}^{\infty }$ (as defined by (%
\ref{E1.1}) and (\ref{E1.2})) converges $\mu $-a.e. on $\Omega $ when we
specialize $f$ and $g$ to be $\mu $-integrable simple functions such that $%
\left\| f\right\| _{L^{\infty }\left( \mu \right) }=\left\| g\right\|
_{L^{\infty }\left( \mu \right) }=1$. This will be carried out in two parts.
(In the ensuing discussion, we shall use without explicit mention the
convenient fact that, since $U\left| L^{\infty }\left( \mu \right) \right. $
is a surjective linear isometry, we have for all $n\in \mathbb{Z}$, 
\begin{equation*}
\left\| U^{n}f\right\| _{L^{\infty }\left( \mu \right) }=\left\|
U^{n}g\right\| _{L^{\infty }\left( \mu \right) }=1\text{.)}
\end{equation*}

Part (i). We first prove the $\mu $-a.e. convergence of $\left\{
A_{k,U}\left( f,g\right) \right\} _{k=1}^{\infty }$ for such $f$ and $g$.
Let $M$ be an arbitrary positive real number, and choose a real-valued
function $K^{\left( M\right) }\in C^{\infty }\left( \mathbb{R}\right) $ such
that: $K^{\left( M\right) }$ vanishes on $\left( -\infty ,-1/M\right]
\bigcup \left[ 1+\dfrac{1}{M},\infty \right) $; $K^{\left( M\right) }=1$ on $%
\left[ 0,1\right] $; $K^{\left( M\right) }$ is increasing on $\left[ -1/M,0%
\right] ;K^{\left( M\right) }$ is decreasing on $\left[ 1,1+\dfrac{1}{M}%
\right] $. By Corollary \ref{Transferred_Oscillation}, for each $m\in 
\mathbb{N}$, the sequence $\left\{ \mathfrak{A}_{j,K^{\left( M\right)
},m,U}\left( f,g\right) \right\} _{j=1}^{\infty }$ corresponding to $f$, $g$%
, $K^{\left( M\right) }$, and $m$ in accordance with (\ref{e2a.32})
converges $\mu $-a.e. on $\Omega $. For each $j\in \mathbb{N}$, we have%
\begin{eqnarray}
&&\hspace{0.15in}  \label{e3.1} \\
&&\,\,\hspace{-0.15in}\hspace{0.15in}\mathfrak{A}_{j,K^{\left( M\right)
},m,U}\left( f,g\right)  \notag \\
&=&\sum \left\{ \left( U^{n}f\right) \left( U^{-n}g\right) K_{j,m}^{\left(
M\right) }\left( n\right) :-\frac{2^{j/m}}{M}<n<2^{j/m}\left( 1+\dfrac{1}{M}%
\right) \right\}  \notag \\
&=&\sum \left\{ \left( U^{n}f\right) \left( U^{-n}g\right) K_{j,m}^{\left(
M\right) }\left( n\right) :-\frac{2^{j/m}}{M}<n<0\right\} +\frac{1}{2^{j/m}}%
\sum_{n=0}^{\left[ 2^{j/m}\right] }\left( U^{n}f\right) \left( U^{-n}g\right)
\notag \\
&&\,\hspace{-0.15in}\hspace{0.15in}\hspace{0.15in}+\sum \left\{ \left(
U^{n}f\right) \left( U^{-n}g\right) K_{j,m}^{\left( M\right) }\left(
n\right) :2^{j/m}<n<2^{j/m}\left( 1+\dfrac{1}{M}\right) \right\} \text{.} 
\notag
\end{eqnarray}%
But for each $n\in \mathbb{Z}$ such that $-\dfrac{2^{j/m}}{M}<n<0$ or $%
2^{j/m}<n<2^{j/m}\left( 1+\dfrac{1}{M}\right) $,%
\begin{equation*}
0\leq K_{j,m}^{\left( M\right) }\left( n\right) \leq \frac{1}{2^{j/m}}\text{,%
}
\end{equation*}%
and so%
\begin{gather}
\left| \sum \left\{ \left( U^{n}f\right) \left( U^{-n}g\right)
K_{j,m}^{\left( M\right) }\left( n\right) :-\frac{2^{j/m}}{M}<n<0\right\}
\right| \leq \frac{1}{2^{j/m}}\left( \frac{2^{j/m}}{M}\right) =\dfrac{1}{M}%
\text{;}  \label{e3.2} \\
\left| \sum \left\{ \left( U^{n}f\right) \left( U^{-n}g\right)
K_{j,m}^{\left( M\right) }\left( n\right) :2^{j/m}<n<2^{j/m}\left( 1+\dfrac{1%
}{M}\right) \right\} \right| \leq \left( \frac{1}{M}+\frac{1}{2^{j/m}}%
\right) \text{.}  \notag
\end{gather}%
Since $\left\{ \mathfrak{A}_{j,K^{\left( M\right) },m,U}\left( f,g\right)
\right\} _{j=1}^{\infty }$ converges $\mu $-a.e. on $\Omega $ for arbitrary $%
m\in \mathbb{N}$, and arbitrary positive $M\in \mathbb{R}$, it follows
readily from (\ref{e3.1}) and (\ref{e3.2}) that for each $m\in \mathbb{N}$
the sequence $\left\{ \widetilde{A}_{2^{j/m},U}\left( f,g\right) \right\}
_{j=1}^{\infty }$ specified by%
\begin{equation*}
\widetilde{A}_{2^{j/m},U}\left( f,g\right) =\frac{1}{2^{j/m}}\sum_{n=0}^{%
\left[ 2^{j/m}\right] -1}\left( U^{n}f\right) \left( U^{-n}g\right) \text{,
for each }j\in \mathbb{N}\text{,}
\end{equation*}%
converges $\mu $-a.e. on $\Omega $. So for Part (i) of the proof it remains
to show that the $\mu $-a.e. convergence of $\left\{ \widetilde{A}%
_{2^{j/m},U}\left( f,g\right) \right\} _{j=1}^{\infty }$ for each fixed $%
m\in \mathbb{N}$ can be converted into $\mu $-a.e. convergence of the
sequence $\left\{ A_{k,U}\left( f,g\right) \right\} _{k=1}^{\infty }$. Given 
$m\in \mathbb{N}$, $k\in \mathbb{N}$, with $k\geq 2$, let $j=j\left(
k,m\right) \in \mathbb{N}$ satisfy%
\begin{equation}
2^{j/m}\leq k<2^{(j+1)/m}\text{.}  \label{e3.11a}
\end{equation}%
Hence for some absolute constant $\eta $ we have 
\begin{equation}
0\leq \frac{k-2^{j/m}}{k}\leq \frac{k-2^{j/m}}{2^{j/m}}<2^{1/m}-1\leq \frac{%
\eta }{m}\text{,}  \label{e3.11b}
\end{equation}%
and consequently we have pointwise on $\Omega $,%
\begin{eqnarray}
&&\left| A_{k,U}\left( f,g\right) -\widetilde{A}_{2^{j/m},U}\left(
f,g\right) \right|  \label{e3.12} \\
&\leq &\left( \frac{k-2^{j/m}}{2^{j/m}k}\right) \sum_{n=0}^{\left[ 2^{j/m}%
\right] -1}\,\left| \left( U^{n}f\right) \right| \,\left| U^{-n}g\right| +%
\frac{1}{k}\sum_{n=\left[ 2^{j/m}\right] }^{k-1}\,\left| \left(
U^{n}f\right) \right| \,\left| U^{-n}g\right|  \notag \\
&\leq &2\left( \frac{k-2^{j/m}}{k}\right) +\frac{2^{j/m}-\left[ 2^{j/m}%
\right] }{k}  \notag \\
&\leq &\frac{2\eta }{m}+\frac{1}{k}\text{.}  \notag
\end{eqnarray}%
It follows from (\ref{e3.12}) and the $\mu $-a.e. convergence for each $m\in 
\mathbb{N}$ of the sequence $\left\{ \widetilde{A}_{2^{j/m},U}\left(
f,g\right) \right\} _{j=1}^{\infty }$ that the sequence $\left\{
A_{k,U}\left( f,g\right) \right\} _{k=1}^{\infty }$ is pointwise Cauchy $\mu 
$-a.e on $\Omega $.

Part (ii). To complete the demonstration of Theorem \ref{Main_Result}, it
will suffice (as noted above) to establish the $\mu $-a.e. convergence of
the averages $\left\{ H_{k,U}\left( f,g\right) \right\} _{k=1}^{\infty }$
for $\mu $-integrable simple functions $f,g$ such that $\left\| f\right\|
_{L^{\infty }\left( \mu \right) }=\left\| g\right\| _{L^{\infty }\left( \mu
\right) }=1$. For this purpose, we shall follow the main outlines of the
proof for Theorem 1.2 of \ \cite{Dem}. We start by letting $M$ be an
arbitrary integer such that $M\geq 2$, and then choosing, as the relevant
kernel for applying Theorem \ref{Dem_Thm_1_4}, an odd $C^{\infty }\left( 
\mathbb{R}\right) $ function $\mathfrak{K}^{\left\langle M\right\rangle }:$ $%
\mathbb{R}\rightarrow \mathbb{R}$ such that:%
\begin{eqnarray}
\mathfrak{K}^{\left\langle M\right\rangle }\left( x\right) &=&\frac{1}{x}%
\text{, for }\left| x\right| \geq 1\text{;}  \label{e3.13} \\
\mathfrak{K}^{\left\langle M\right\rangle }\left( x\right) &=&0\text{, for }%
\left| x\right| \leq 1-\frac{1}{M}\text{;}  \label{e3.13b} \\
\left| \mathfrak{K}^{\left\langle M\right\rangle }\left( x\right) \right|
&\leq &2\text{, for }\left| x\right| \leq 1\text{.}  \label{e3.13c}
\end{eqnarray}%
Notice that by virtue of (\ref{e3.13}) , the successive derivatives $\dfrac{%
d^{n}\mathfrak{K}^{\left\langle M\right\rangle }\left( x\right) }{dx^{n}}$,
for $n\in \mathbb{N}$, all belong to $L^{1}\left( \mathbb{R}\right) \bigcap
L^{\infty }\left( \mathbb{R}\right) $, and this fact is helpful in seeing by
elementary considerations that $\mathfrak{K}^{\left\langle M\right\rangle }$
satisfies the hypotheses of Theorem \ref{Dem_Thm_1_4}, which thereby
furnishes us with the following oscillation estimate, valid for each integer 
$M\geq 2$, each $m\in \mathbb{N}$, each integer $J\geq 2$, each sequence of
positive integers $u_{1}<u_{2}<\cdots <u_{J}$, and every pair of compactly
supported functions $F$ and $G$ belonging to $L^{\infty }\left( \mathbb{R}%
\right) $.%
\begin{align}
& \left\| \left\{ \sum_{j=1}^{J-1}\sup_{\substack{ k\in \mathbb{N},  \\ %
u_{j}\leq k<u_{j+1}}}\left| \int_{\mathbb{R}}F\left( x+y\right) G\left(
x-y\right) \left( \left( \delta _{d_{m}^{k}}\mathfrak{K}^{\left\langle
M\right\rangle }\right) \left( y\right) -\right. \right. \right. \right.
\label{e3.14} \\
& \hspace{0.5in}\hspace{0.5in}\hspace{0.5in}\hspace{0.5in}\left. \left.
\left. \left. \left( \delta _{d_{m}^{u_{j+1}}}\mathfrak{K}^{\left\langle
M\right\rangle }\right) \left( y\right) \right) \,dy\right| ^{2}\right\}
^{1/2}\right\| _{L_{x}^{1,\infty }}  \notag \\
& \leq C_{M,m}\,J^{1/4}\,\left\| F\right\| _{L^{2}\left( \mathbb{R}\right)
}\,\left\| G\right\| _{L^{2}\left( \mathbb{R}\right) }\text{.}  \notag
\end{align}%
However, since the kernel $\mathfrak{K}^{\left\langle M\right\rangle }$
lacks compact support and does not belong to $L^{1}\left( \mathbb{R}\right) $%
, its exploitation of the oscillation estimate (\ref{e3.14}) will require
extra care. In this regard, it is convenient to observe that because of (\ref%
{e3.13}) $\mathfrak{K}^{\left\langle M\right\rangle }$ has the following
``quasi-stability'' under dilations: for each positive real number $\xi $, 
\begin{equation*}
\left( \delta _{\xi }\mathfrak{K}^{\left\langle M\right\rangle }\right)
\left( x\right) =\frac{1}{x}\text{, whenever }\left| x\right| \geq \xi \text{%
.}
\end{equation*}%
Hence if $0<\xi _{1}\leq $ $\xi _{2}$, then $\left( \delta _{\xi _{1}}%
\mathfrak{K}^{\left\langle M\right\rangle }\right) \left( x\right) =\left(
\delta _{\xi _{2}}\mathfrak{K}^{\left\langle M\right\rangle }\right) \left(
x\right) =\dfrac{1}{x}$, for $\left| x\right| \geq \xi _{2}$. In particular,
for $1\leq j\leq J-1$, and $u_{j}\leq k<u_{j+1}$, 
\begin{equation}
\left( \delta _{d_{m}^{k}}\mathfrak{K}^{\left\langle M\right\rangle }-\delta
_{d_{m}^{u_{j+1}}}\mathfrak{K}^{\left\langle M\right\rangle }\right) \left(
x\right) =0\text{, whenever }\left| x\right| >d_{m}^{u_{j+1}}\text{.}
\label{e3.14a}
\end{equation}%
Consequently, although the $C^{\infty }\left( \mathbb{R}\right) $ kernel $%
\mathfrak{K}^{\left\langle M\right\rangle }$ itself lacks compact support
(and so is not automatically covered by the discretization result in Theorem %
\ref{discrete_Thm1_4}), the discretization methods in the proof of Theorem %
\ref{discrete_Thm1_4} can nevertheless go forward from (\ref{e3.14}) by
straightforward adjustments which rely on the compact supports of the
relevant difference kernels in accordance with (\ref{e3.14a}). This
procedure discretizes (\ref{e3.14}) by yielding the following result. For
each integer $M\geq 2$, for each $m\in \mathbb{N}$, for every pair of
finitely supported sequences $a\in \mathbb{C}^{\mathbb{Z}}$ and $b\in 
\mathbb{C}^{\mathbb{Z}}$, for every integer $R\geq 2$, and for each sequence
of positive integers $u_{1}<u_{2}<\cdots <u_{R}$, 
\begin{eqnarray}
&&\left\| \left\{ \sum_{r=1}^{R-1}\sup_{\substack{ j\in \mathbb{N},  \\ %
u_{r}\leq j<u_{r+1}}}\left| \mathcal{K}_{j,M,m}\left( a,b\right) -\mathcal{K}%
_{u_{r+1},M,m}\left( a,b\right) \right| ^{2}\right\} ^{%
{\frac12}%
}\right\| _{\ell ^{1,\infty }\left( \mathbb{Z}\right) }  \label{e3.15} \\
&\leq &C_{M,m}\,R^{1/4}\,\left\| a\right\| _{\ell ^{2}\left( \mathbb{Z}%
\right) }\left\| b\right\| _{\ell ^{2}\left( \mathbb{Z}\right) }\text{,} 
\notag
\end{eqnarray}%
where for each $\nu \in \mathbb{N}$, and each $k\in \mathbb{Z}$,%
\begin{equation}
\left( \mathcal{K}_{\nu ,M,m}\left( a,b\right) \right) \left( k\right)
=\sum_{n=-\infty }^{\infty }\,\frac{a\left( k+n\right) \,b\left( k-n\right) 
}{d_{m}^{\nu }}\mathfrak{K}^{\left\langle M\right\rangle }\left( \frac{n}{%
d_{m}^{\nu }}\right) .  \label{e3.16}
\end{equation}%
(Since $a$ and $b$ are finitely supported, the sum on the right of (\ref%
{e3.16}) has only finitely many non-zero terms, and also $\mathcal{K}_{\nu
,M,m}\left( a,b\right) $ is finitely supported.)

In order to obtain a transferred counterpart of (\ref{e3.15}) to which we
can apply Theorem \ref{gen_discrete_tran}, we shall first recast (\ref{e3.15}%
) so that it becomes completely expressed in terms of finitely supported
discrete kernels (rather than the present discrete kernels of the form $%
\left( \delta _{d_{m}^{\nu }}\mathfrak{K}^{\left\langle M\right\rangle
}\right) \left| \mathbb{Z}\right. $). The method for doing so will be taken
from the proof \ for Theorem 1.2 in \cite{Dem}. Specifically, for each $j\in 
\mathbb{N}$, we define $\mathbb{A}_{j,M,m}:$ $\mathbb{Z}\rightarrow \mathbb{R%
}$, $\mathfrak{H}_{j,m}:\mathbb{Z}\rightarrow \mathbb{R}$ and $\mathfrak{D}%
_{j,M,m}:\mathbb{Z}\rightarrow \mathbb{R}$ by writing:%
\begin{eqnarray}
\mathbb{A}_{j,M,m}\left( n\right) &=&\left\{ 
\begin{array}{cc}
\frac{1}{d_{m}^{j}}\mathfrak{K}^{\left\langle M\right\rangle }\left( \frac{n%
}{d_{m}^{j}}\right) \text{,} & \text{if }\left| n\right| \leq d_{m}^{j}\text{%
;} \\ 
0\text{,} & \text{otherwise.}%
\end{array}%
\right.  \label{e3.16a} \\
\mathfrak{H}_{j,m}\left( n\right) &=&\left\{ 
\begin{array}{cc}
\frac{1}{n}\text{,} & \text{if }0<\left| n\right| \leq d_{m}^{j}\text{;} \\ 
0\text{,} & \text{otherwise.}%
\end{array}%
\right.  \notag \\
\mathfrak{D}_{j,M,m}\left( n\right) &=&\mathbb{A}_{j,M,m}\left( n\right) -%
\mathfrak{H}_{j,m}\left( n\right) \text{, for all }n\in \mathbb{Z}\text{.} 
\notag
\end{eqnarray}%
Then it is easy to verify from definitions that whenever $j\in \mathbb{N}$
and $\nu \in \mathbb{N}$, we have for all $n\in \mathbb{Z}$, 
\begin{eqnarray*}
&&\frac{1}{d_{m}^{j}}\mathfrak{K}^{\left\langle M\right\rangle }\left( \frac{%
n}{d_{m}^{j}}\right) -\frac{1}{d_{m}^{\nu }}\mathfrak{K}^{\left\langle
M\right\rangle }\left( \frac{n}{d_{m}^{\nu }}\right) \\
&=&\mathfrak{D}_{j,M,m}\left( n\right) -\mathfrak{D}_{\nu ,M,m}\left(
n\right) \text{.}
\end{eqnarray*}%
For each $\nu \in \mathbb{N}$, we define the bilinear mapping $D_{\nu ,M,m}:%
\mathbb{C}^{\mathbb{Z}}\,\times \,\mathbb{C}^{\mathbb{Z}}\rightarrow \mathbb{%
C}^{\mathbb{Z}}$ by writing for all $\mathfrak{v}\in \mathbb{C}^{\mathbb{Z}}$%
, all $\mathfrak{w}\in \mathbb{C}^{\mathbb{Z}}$, and all $k\in \mathbb{Z}$, 
\begin{equation*}
\left( D_{\nu ,M,m}\left( \mathfrak{v},\mathfrak{w}\right) \right) \left(
k\right) =\sum_{n=-\infty }^{\infty }\mathfrak{v}\left( k+n\right) \,%
\mathfrak{w}\left( k-n\right) \,\mathfrak{D}_{\nu ,M,m}\left( n\right) \text{%
.}
\end{equation*}%
We can now use the finitely supported discrete kernels $\mathfrak{D}_{j,M,m}$
to rewrite the inequality (\ref{e3.15}) in the following form, valid for
each integer $M\geq 2$, each $m\in \mathbb{N}$, every pair of finitely
supported sequences $a\in \mathbb{C}^{\mathbb{Z}}$ and $b\in \mathbb{C}^{%
\mathbb{Z}}$, every integer $R\geq 2$, and each sequence of positive
integers $u_{1}<u_{2}<\cdots <u_{R}$.%
\begin{eqnarray}
&&\left\| \left\{ \sum_{r=1}^{R-1}\sup_{\substack{ j\in \mathbb{N},  \\ %
u_{r}\leq j<u_{r+1}}}\left| D_{j,M,m}\left( a,b\right)
-D_{u_{r+1},M,m}\left( a,b\right) \right| ^{2}\right\} ^{%
{\frac12}%
}\right\| _{\ell ^{1,\infty }\left( \mathbb{Z}\right) }  \label{e3.17} \\
&\leq &C_{M,m}\,R^{1/4}\,\left\| a\right\| _{\ell ^{2}\left( \mathbb{Z}%
\right) }\left\| b\right\| _{\ell ^{2}\left( \mathbb{Z}\right) }\text{.} 
\notag
\end{eqnarray}%
After applying Theorem \ref{gen_discrete_tran} to (\ref{e3.17}), we infer
that for each integer $M\geq 2$, for each $m\in \mathbb{N}$, and for the
above-described $\mu $-integrable simple functions $f,g$, the sequence%
\begin{equation}
\left\{ \sum_{\left| k\right| \leq d_{m}^{n}}\left( U^{k}f\right) \,\left(
U^{-k}g\right) \mathfrak{D}_{n,M,m}\left( k\right) \right\} _{n=1}^{\infty }%
\text{ converges }\mu \text{-a.e. on }\Omega \text{.}  \label{e3.18}
\end{equation}%
For each $n\in \mathbb{N}$, it is clear from definitions and the notation of
(\ref{E1.2}) that the following identity holds pointwise on $\Omega $.%
\begin{equation}
\sum_{\left| k\right| \leq d_{m}^{n}}\left( U^{k}f\right) \,\left(
U^{-k}g\right) \mathfrak{D}_{n,M,m}\left( k\right) =\sum_{\left| k\right|
\leq d_{m}^{n}}\left( U^{k}f\right) \,\left( U^{-k}g\right) \mathbb{A}%
_{n,M,m}\left( k\right) -H_{d_{m}^{n},U}\left( f,g\right) \text{.}
\label{e3.19}
\end{equation}%
Moreover, (\ref{e3.16a}), taken in conjunction with (\ref{e3.13b}) and (\ref%
{e3.13c}), shows that for each integer $M\geq 2$, each $m\in \mathbb{N}$,
and each $n\in \mathbb{N}$, 
\begin{eqnarray*}
&&\sum_{\left| k\right| \leq d_{m}^{n}}\left( U^{k}f\right) \,\left(
U^{-k}g\right) \mathbb{A}_{n,M,m}\left( k\right) \\
&=&\sum \left\{ \left( U^{k}f\right) \,\left( U^{-k}g\right) \frac{1}{%
d_{m}^{n}}\mathfrak{K}^{\left\langle M\right\rangle }\left( \frac{k}{%
d_{m}^{n}}\right) :\left( 1-\frac{1}{M}\right) d_{m}^{n}<\left| k\right|
\leq d_{m}^{n}\right\} \text{,}
\end{eqnarray*}%
and so we have pointwise on $\Omega $, 
\begin{equation*}
\left| \sum_{\left| k\right| \leq d_{m}^{n}}\left( U^{k}f\right) \,\left(
U^{-k}g\right) \mathbb{A}_{n,M,m}\left( k\right) \right| \leq 4\left( \frac{1%
}{M}+\frac{1}{d_{m}^{n}}\right) =4\left( \frac{1}{M}+\frac{1}{2^{n/m}}%
\right) \text{.}
\end{equation*}%
Since the integer $M\geq 2$ is arbitrary, it follows from this, (\ref{e3.18}%
), and (\ref{e3.19}) that for each $m\in \mathbb{N}$,%
\begin{equation}
\left\{ H_{2^{n/m},U}\left( f,g\right) \right\} _{n=1}^{\infty }\text{
converges }\mu \text{-a.e. on }\Omega \text{.}  \label{e3.20}
\end{equation}

\noindent The $\mu $-a.e. convergence of the averages $\left\{ H_{k,U}\left(
f,g\right) \right\} _{k=1}^{\infty }$ can now be deduced from (\ref{e3.20})
by reasoning analogous to that used to complete Part (i) of the proof. \ \ $%
\square $

\section{The Continuous Variable Counterpart of Theorem \ref{Main_Result}%
\label{sec4}}

This brief final section features Theorem \ref{one-parameter-ver} below,
which is the counterpart of Theorem \ref{Main_Result} for averages defined
by one-parameter groups of Lebesgue space isometries associated with the
arbitrary sigma-finite measure space $\left( \Omega ,\mu \right) $. \ These
continuous variable averages do not require any discretization for
oscillation estimates based on Theorem \ref{Dem_Thm_1_4}, and in this
respect their treatment is simpler than that for the discrete averages.
Moreover, Theorem \ref{one-parameter-ver} below can be established by
techniques which, though at times involving measure-theoretic
technicalities, are transparently analogous to those used in the preceding
sections for the discrete averages. In particular, for the relevant
transferred maximal estimates (repectively, relevant transferred oscillation
estimates) in the present setting, we need only replace the role of Theorem %
\ref{bisublinear_maxml} (respectively, Theorem \ref{gen_discrete_tran}) by
suitable reasoning based on \S 6 of \cite{BBCG} so as to transfer from $%
\mathbb{R}$ to the $\left( \Omega ,\mu \right) $ context Michael Lacey's
classical estimates (\cite{Lacey}) for the bisublinear Hardy-Littlewood
maximal operator and the bisublinear maximal Hilbert transform
(respectively, oscillation estimates of the form (\ref{e2a.0}) of Theorem %
\ref{Dem_Thm_1_4}). In view of this state of affairs, the discussion below
will, for expository reasons, omit detailed arguments.

In order to formulate the results for the continuous variable setting, we
begin by describing the main ingredients of the discussion. Our transference
vehicle for defining the relevant averages on the measure space side will be
a one-parameter group $\mathcal{U}\equiv \left\{ U_{t}:t\in \mathbb{R}%
\right\} $ consisting of linear bijections of $\mathcal{A}\left( \mu \right) 
$ onto $\mathcal{A}\left( \mu \right) $. Thus,%
\begin{equation}
U_{s+t}\left( f\right) =U_{s}\left( U_{t}f\right) \text{, for all }s\in 
\mathbb{R}\text{, }t\in \mathbb{R}\text{, }f\in \mathcal{A}\left( \mu
\right) \text{.}  \label{e4.a}
\end{equation}%
The one-parameter group $\mathcal{U}\equiv \left\{ U_{t}:t\in \mathbb{R}%
\right\} $ will be required to satisfy the following conditions.

\begin{enumerate}
\item[(C1)] For each $t\in \mathbb{R}$, $\lim_{k\rightarrow \infty }\,\left(
U_{t}g_{k}\right) =U_{t}g$ $\mu $-a.e. on $\Omega $, whenever $\left\{
g_{k}\right\} _{k=1}^{\infty }\subseteq \mathcal{A}\left( \mu \right) $, $%
g\in \mathcal{A}\left( \mu \right) $, and $\lim_{k\rightarrow \infty
}\,g_{k}=g$ $\mu $-a.e. on $\Omega $.

\item[(C2)] For $0<p<\infty $, and each $s\in \mathbb{R}$, $L^{p}\left( \mu
\right) $ is invariant under $U_{s}$, and the restrictions $\left\{
U_{t}\left| L^{p}\left( \mu \right) \right. :t\in \mathbb{R}\right\} $ form
a strongly continuous one-parameter group of surjective linear isometries of 
$L^{p}\left( \mu \right) $ onto $L^{p}\left( \mu \right) $.

\item[(C3)] For each $f\in \mathcal{A}\left( \mu \right) $, the expression $%
\left( U_{t}f\right) \left( x\right) $, where $\left( t,x\right) $ runs
through $\mathbb{R}\,\times \,\Omega $, can be regarded as being a jointly
measurable version with respect to the product of linear Lebesgue measure $%
\mathfrak{m}_{\mathbb{R}}$ and the measure $\mu $. In other words, there
exists a complex-valued $\left( \mathfrak{m}_{\mathbb{R}}\,\times \mu
\right) $-measurable function $\mathfrak{F}_{f}$ on $\mathbb{R}\,\times
\,\Omega $ such that for each $t\in \mathbb{R}$, $\mathfrak{F}_{f}\left(
t,\bullet \right) $ belongs to the equivalence class (modulo equality $\mu $%
-a.e. on $\Omega $) of $U_{t}f$. (For convenience, we shall denote such a
function $\mathfrak{F}_{f}$ by $\left( U_{t}f\right) \left( x\right) $.)
\end{enumerate}

\begin{remark}
\label{re_props_of_U_t}(i)We observe here that for $f\in \mathcal{A}\left(
\mu \right) $, any two jointly measurable versions $\mathfrak{F}_{f}^{\left(
1\right) }$ \ and $\mathfrak{F}_{f}^{\left( 2\right) }$ representing $\left(
U_{t}f\right) \left( x\right) $ on $\mathbb{R}\,\times \,\Omega $ \ in
accordance with condition (C3) would automatically have the additional
property that for $\mu $-almost all $x\in \Omega $,%
\begin{equation}
\mathfrak{F}_{f}^{\left( 1\right) }\left( t,x\right) =\mathfrak{F}%
_{f}^{\left( 2\right) }\left( t,x\right) \text{, for }\mathfrak{m}_{\mathbb{R%
}}\text{-almost all }t\in \mathbb{R}\text{.}  \label{e4.aaa}
\end{equation}%
For this reason among other obvious reasons, the particular choice of
jointly measurable version of $\left( U_{t}f\right) \left( x\right) $ on $%
\mathbb{R}\,\times \,\Omega $ will be immaterial in all our considerations
below.

(ii)As is well-known (see, e.g., Proposition 5 in \cite{BBCG}), the above
conditions (C1) and (C2) can be shown to imply that for each $t\in \mathbb{R}$%
, $U_{t}\left| L^{\infty }\left( \mu \right) \right. $ is a surjective
linear isometry of $L^{\infty }\left( \mu \right) $ onto $L^{\infty }\left(
\mu \right) $. Hence, in view of condition (C3), for every $f\in $ $L^{\infty
}\left( \mu \right) $, we have for $\mu $-almost all $x\in \Omega $,%
\begin{equation}
\left\| \left( U_{\left( \cdot \right) }f\right) \left( x\right) \right\|
_{L^{\infty }\left( \mathbb{R}\right) }\leq \left\| f\right\| _{L^{\infty
}\left( \mu \right) }\text{.}  \label{e4aa}
\end{equation}
\end{remark}

In terms of the foregoing notation, our continuous variable version of
Theorem \ref{Main_Result} takes the following form. (In particular, the
variable limits of integration of the indefinite integrals occurring below
are continuous rather than discrete variables.)

\begin{theorem}
\label{one-parameter-ver}Let $\left( \Omega ,\mu \right) $ be a sigma-finite
measure space, let $\mathcal{U}\equiv \left\{ U_{t}:t\in \mathbb{R}\right\} $
be a one-parameter group of linear bijections of $\mathcal{A}\left( \mu
\right) $ onto $\mathcal{A}\left( \mu \right) $ satisfying the above
conditions (C1), (C2), and (C3), and let $p_{1}$, $p_{2}$, $p_{3}$ satisfy (\ref%
{E1.3}) and (\ref{E1.4}). Then for every pair of functions $f\in
L^{p_{1}}\left( \mu \right) \bigcap $ $L^{2}\left( \mu \right) $ and $g\in
L^{p_{2}}\left( \mu \right) \bigcap $ $L^{2}\left( \mu \right) $, the
following assertions hold.

\begin{enumerate}
\item[(i)] Each of the following two limits exists $\mu $-a.e. on $\Omega $,
as well as with respect to the metric topology of $\ $the space $%
L^{p_{3}}\left( \mu \right) $. 
\begin{eqnarray}
\left( \mathbb{E}_{\mathcal{U},\infty }\left( f,g\right) \right) \left(
x\right) &\equiv &\lim_{r\rightarrow \infty }\,\frac{1}{r}\int_{0}^{r}\left(
U_{t}f\right) \left( x\right) \left( U_{-t}g\right) \left( x\right) \,dt%
\text{;}  \label{e4.11a} \\
\left( \mathbb{E}_{\mathcal{U},0}\left( f,g\right) \right) \left( x\right)
&\equiv &\lim_{r\rightarrow 0^{+}}\,\frac{1}{r}\int_{0}^{r}\left(
U_{t}f\right) \left( x\right) \left( U_{-t}g\right) \left( x\right) \,dt%
\text{.}  \label{e4.11b}
\end{eqnarray}%
Moreover,%
\begin{equation}
\max \left\{ \left\| \mathbb{E}_{\mathcal{U},0}\left( f,g\right) \right\|
_{L^{p_{3}}\left( \mu \right) }\text{,}\left\| \mathbb{E}_{\mathcal{U}%
,\infty }\left( f,g\right) \right\| _{L^{p_{3}}\left( \mu \right) }\right\}
\leq C_{p_{1},p_{2}}\left\| f\right\| _{L^{p_{1}}\left( \mu \right)
}\,\left\| g\right\| _{L^{p_{2}}\left( \mu \right) }\text{.}  \label{e4.12}
\end{equation}

\item[(ii)] For $\mu $-almost all $x\in \Omega $, the (Cauchy principal
value) improper integral%
\begin{equation}
\int_{\varepsilon \leq \left| t\right| }\dfrac{\left( U_{t}f\right) \left(
x\right) \left( U_{-t}g\right) \left( x\right) }{t}\,dt\equiv
\lim_{\varsigma \rightarrow \infty }\int_{\varepsilon \leq \left| t\right|
\leq \varsigma }\dfrac{\left( U_{t}f\right) \left( x\right) \left(
U_{-t}g\right) \left( x\right) }{t}\,dt  \label{e4.12aa}
\end{equation}%
exists in $\mathbb{C}$ for every $\varepsilon >0$, and $\left( \mathbb{H}%
_{\varepsilon ,\mathcal{U}}\left( f,g\right) \right) \left( x\right) \equiv
\int_{\varepsilon \leq \left| t\right| }\dfrac{\left( U_{t}f\right) \left(
x\right) \left( U_{-t}g\right) \left( x\right) }{t}\,dt$ approaches a limit $%
\left( \mathbb{H}_{\mathcal{U}}\left( f,g\right) \right) \left( x\right) \in 
\mathbb{C}$, as $\varepsilon \rightarrow 0^{+}$. We also have the following
two limit relations with respect to convergence in the metric topology of $\ 
$the space $L^{p_{3}}\left( \mu \right) $.%
\begin{gather}
\mathbb{H}_{\varepsilon ,\mathcal{U}}\left( f,g\right) =\lim_{\varsigma
\rightarrow \infty }\,\int_{\varepsilon \leq \left| t\right| \leq \varsigma }%
\dfrac{\left( U_{t}f\right) \left( \cdot \right) \left( U_{-t}g\right)
\left( \cdot \right) }{t}\,dt\text{, }  \label{e4.12a} \\
\hspace{0.8in}\hspace{0.8in}\hspace{0.5in}\hspace{0.8in}\hspace{0.5in}\text{%
for each }\varepsilon >0\text{;}  \notag \\
\mathbb{H}_{\mathcal{U}}\left( f,g\right) =\lim_{\varepsilon \rightarrow
0^{+}}\,\mathbb{H}_{\varepsilon ,\mathcal{U}}\left( f,g\right) \text{.}
\label{e4.13}
\end{gather}%
Moreover, 
\begin{equation}
\left\| \mathbb{H}_{\mathcal{U}}\left( f,g\right) \right\| _{L^{p_{3}}\left(
\mu \right) }\leq \left\| \sup_{\varepsilon >0}\,\left| \mathbb{H}%
_{\varepsilon ,\mathcal{U}}\left( f,g\right) \right| \right\|
_{L^{p_{3}}\left( \mu \right) }\leq C_{p_{1},p_{2}}\left\| f\right\|
_{L^{p_{1}}\left( \mu \right) }\,\left\| g\right\| _{L^{p_{2}}\left( \mu
\right) }\text{. }  \label{e4.13a}
\end{equation}
\end{enumerate}
\end{theorem}

We come now to a discussion of the convergence properties exhibited by the
one-parameter (continuous variable) averages in the setting of an arbitrary
measure space $\left( X,\sigma \right) $. Although the measure $\sigma $
need not be sigma-finite, this context does offer a notion of joint
measurability with respect to the measurable spaces of $\mathfrak{m}_{%
\mathbb{R}}$ and $\sigma $ for complex-valued functions defined on $\mathbb{R%
}\,\times X\,$(as in, e.g., \S 33 of \cite{Halmos}), but a product measure
of $\mathfrak{m}_{\mathbb{R}}$ and $\sigma $ (in the sense of the abstract
Fubini's theorem, as in, e.g., \S \S 35,36 of \cite{Halmos}) is lacking, and
this lack imposes technical constraints on attempts to mirror the general
measure space results for the a.e. convergence of the discrete averages
induced by measure-preserving point transformations (Corollary \ref%
{4all_measure_spaces}). For example, we no longer have a route to
compatibility conditions like\ (\ref{e4.aaa}), and so the framing of $\sigma 
$-a.e. convergence questions in the continuous variable framework can become
refractory. In view of these circumstances we shall, for convenience, forgo
discussion of $\sigma $-a.e. convergence for the one-parameter averages in
favor of studying the convergence in $L^{p}\left( \sigma \right) $ of their
Bochner integral formulations.

The transference vehicle for the present framework will be a one-parameter
group (under composition of mappings) $\mathcal{P}\equiv \left\{ \psi
_{t}:t\in \mathbb{R}\right\} $ consisting of invertible measure-preserving
point transformations of $\left( X,\sigma \right) $--thus, for all $x\in X$,
all $s\in \mathbb{R}$, and all $t\in \mathbb{R}$, $\psi _{s+t}\left(
x\right) =\psi _{s}\left( \psi _{t}\left( x\right) \right) $. In this setup,
we postulate the following two properties for $\mathcal{P}\equiv \left\{
\psi _{t}:t\in \mathbb{R}\right\} $, thereby endowing $\mathcal{V}\equiv
\left\{ V_{t}:t\in \mathbb{R}\right\} $, the one-parameter group of
corresponding composition operators on $\mathcal{A}\left( \sigma \right) $,
with the counterpart of the conditions (C1), (C2), and (C3) that were imposed
on $\mathcal{U\equiv }\left\{ U_{t}:t\in \mathbb{R}\right\} $ at the outset
of this section.

\begin{enumerate}
\item[(A)] $\lim_{t\rightarrow t_{0}}\sigma \left( \left( \psi _{t}\left(
E\right) \right) \Delta \left( \psi _{t_{0}}\left( E\right) \right) \right)
=0$, for each set $E\subseteq X$ such that $\sigma \left( E\right) <\infty $%
, and each $t_{0}\in \mathbb{R}$,where $\Delta $ denotes the symmetric
difference of sets.

\item[(B)] For each $f\in \mathcal{A}\left( \sigma \right) $, there exists a
complex-valued function $\mathfrak{E}_{f}$ on $\mathbb{R}\,\times \,X$ such
that $\mathfrak{E}_{f}$ is jointly measurable with respect to the measurable
spaces of $\mathfrak{m}_{\mathbb{R}}$ and $\sigma $, and for each $t\in 
\mathbb{R}$, $\mathfrak{E}_{f}\left( t,\bullet \right) $ belongs to the
equivalence class (modulo equality $\sigma $-a.e. on $X$) of $V_{t}f$. (For
convenience, we shall denote such a function $\mathfrak{E}_{f}$ by $\left(
V_{t}f\right) \left( x\right) \equiv f\left( \psi _{t}\left( x\right)
\right) $.)
\end{enumerate}

In view of (A) and Cauchy-Schwarz, for each $f\in L^{2}\left( \sigma \right) 
$ and each $g\in $ $L^{2}\left( \sigma \right) $, the pointwise product $%
f\left( \psi _{t}\right) g\left( \psi _{-t}\right) $ \textit{qua} function
of $t\in \mathbb{R}$ moves continuously in $L^{1}\left( \sigma \right) $. So
for each $F\in L^{1}\left( \mathbb{R}\right) $, the $L^{1}\left( \sigma
\right) $-valued Bochner integral $\int_{\mathbb{R}}\,f\left( \psi
_{t}\right) g\left( \psi _{-t}\right) F\left( t\right) \,dt$ exists and
clearly satisfies%
\begin{equation*}
\left\| \int_{\mathbb{R}}\,f\left( \psi _{t}\right) g\left( \psi
_{-t}\right) F\left( t\right) \,dt\right\| _{L^{1}\left( \sigma \right)
}\leq \left\| F\right\| _{L^{1}\left( \mathbb{R}\right) }\left\| f\right\|
_{L^{2}\left( \sigma \right) }\left\| g\right\| _{L^{2}\left( \sigma \right)
}\text{.}
\end{equation*}

With due attention to technical details arising in this context, one can
deduce, as a corollary of Theorem \ref{one-parameter-ver}, the following
continuous variable variant of Corollary \ref{4all_measure_spaces}.

\begin{corollary}
\label{1param4all_meas-Cor}Suppose that $\left( X,\sigma \right) $ is an
arbitrary measure space, and let $\mathcal{P}\equiv \left\{ \psi _{t}:t\in 
\mathbb{R}\right\} $ be a one-parameter group of invertible
measure-preserving point transformations of $\left( X,\sigma \right) $ onto $%
\left( X,\sigma \right) $ which has the properties (A) and (B) listed above.
Let $p_{1}$, $p_{2}$, $p_{3}$ satisfy (\ref{E1.3}) and (\ref{E1.4}). Then
for each pair of functions $f\in L^{p_{1}}\left( \sigma \right) \bigcap $ $%
L^{2}\left( \sigma \right) $ and $g\in L^{p_{2}}\left( \sigma \right)
\bigcap $ $L^{2}\left( \sigma \right) $, the following assertions are valid.

\begin{enumerate}
\item[(a)] The $L^{1}\left( \sigma \right) $-valued Bochner integrals $%
\int_{0}^{r}f\left( \psi _{t}\right) g\left( \psi _{-t}\right) \,dt$ ($r>0$)
belong to $L^{p_{3}}\left( \sigma \right) $, and have the property that both
the following limits exist with respect to the metric topology of the space $%
L^{p_{3}}\left( \sigma \right) $. 
\begin{eqnarray}
\mathbb{E}_{\mathcal{P},\infty }\left( f,g\right) &\equiv
&\lim_{r\rightarrow \infty }\,\frac{1}{r}\int_{0}^{r}f\left( \psi
_{t}\right) g\left( \psi _{-t}\right) \,dt\text{;}  \label{e4.21a} \\
\mathbb{E}_{\mathcal{P},0}\left( f,g\right) &\equiv &\lim_{r\rightarrow
0^{+}}\,\frac{1}{r}\int_{0}^{r}f\left( \psi _{t}\right) g\left( \psi
_{-t}\right) \,dt\text{.}  \label{e4.21b}
\end{eqnarray}%
Moreover,%
\begin{equation}
\max \left\{ \left\| \mathbb{E}_{\mathcal{P},0}\left( f,g\right) \right\|
_{L^{p_{3}}\left( \sigma \right) },\left\| \mathbb{E}_{\mathcal{P},\infty
}\left( f,g\right) \right\| _{L^{p_{3}}\left( \sigma \right) }\right\} \leq
C_{p_{1},p_{2}}\left\| f\right\| _{L^{p_{1}}\left( \sigma \right) }\,\left\|
g\right\| _{L^{p_{2}}\left( \sigma \right) }\text{.}  \label{e4.21c}
\end{equation}

\item[(b)] For each $\varepsilon >0$, the $L^{1}\left( \sigma \right) $%
-valued Bochner integrals $\int_{\varepsilon \leq \left| t\right| \leq
\varsigma }\,f\left( \psi _{t}\right) g\left( \psi _{-t}\right) t^{-1}\,dt$ (%
$\varepsilon <\varsigma $) belong to $L^{p_{3}}\left( \sigma \right) $, and
have the property that, with respect to the metric topology of the space $%
L^{p_{3}}\left( \sigma \right) $, 
\begin{equation*}
\mathbb{H}_{\varepsilon ,\mathcal{P}}\left( f,g\right) \equiv
\lim_{\varsigma \rightarrow \infty }\int_{\varepsilon \leq \left| t\right|
\leq \varsigma }f\left( \psi _{t}\right) g\left( \psi _{-t}\right) t^{-1}\,dt
\end{equation*}%
exists. We also have, with respect to the metric topology of the space $%
L^{p_{3}}\left( \sigma \right) $, the existence of%
\begin{equation*}
\mathbb{H}_{\mathcal{P}}\left( f,g\right) \equiv \lim_{\varepsilon
\rightarrow 0^{+}}\,\mathbb{H}_{\varepsilon ,\mathcal{P}}\left( f,g\right) 
\text{.}
\end{equation*}%
Moreover,%
\begin{equation*}
\left\| \mathbb{H}_{\mathcal{P}}\left( f,g\right) \right\| _{L^{p_{3}}\left(
\sigma \right) }\leq \sup_{\varepsilon >0}\left\| \mathbb{H}_{\varepsilon ,%
\mathcal{P}}\left( f,g\right) \right\| _{L^{p_{3}}\left( \sigma \right)
}\leq C_{p_{1},p_{2}}\left\| f\right\| _{L^{p_{1}}\left( \sigma \right)
}\,\left\| g\right\| _{L^{p_{2}}\left( \sigma \right) }\text{.}
\end{equation*}
\end{enumerate}
\end{corollary}

We close with the following example, which illustrates Theorem \ref%
{one-parameter-ver} in the realm of harmonic analysis on groups.

\begin{example}
\label{aaexample4Helson-theory} The context of this example will be Helson's
classic theory of generalized analyticity and invariant subspaces, which we
first describe in order to set the stage. (For a full discussion of this
context and its generalizations, we refer the reader to \cite{ABG},\cite%
{Helson}.) Let $\Gamma $ be a dense subgroup of the additive group $\mathbb{R%
}$ of all real numbers. Endow $\Gamma $ with the discrete topology and the
order it inherits from $\mathbb{R}$, and let $K$ be the dual group of $%
\Gamma $. (Equivalently, $K$ can be characterized as a compact abelian group
other than $\left\{ 0\right\} $ or the unit circle $\mathbb{T}$ such that
the dual group of $K$ is archimedean ordered.) \ Denote the normalized Haar
measure of $K$ by $\mathfrak{m}_{K}$, and for each $t\in \mathbb{R}$, let $%
\mathfrak{e}_{t}\in K$ be specified by writing for all $\gamma \in \Gamma $, 
$\mathfrak{e}_{t}\left( \gamma \right) =e^{it\gamma }$. A cocycle on $K$ is
a Borel measurable function $A:\mathbb{R}$\thinspace $\times \,K\rightarrow 
\mathbb{T}$ such that%
\begin{equation*}
A\left( t+u,x\right) =A\left( t,x\right) A\left( u,x+\mathfrak{e}_{t}\right) 
\text{, for all }t\in \mathbb{R}\text{, all }u\in \mathbb{R}\text{, and all }%
x\in K\text{.}
\end{equation*}

Every cocycle $A$ on $K$ automatically has the property that the mapping $%
t\rightarrow A\left( t,\cdot \right) $ is continuous from $\mathbb{R}$ into $%
L^{p}\left( \mathfrak{m}_{K}\right) $ for $0<p<\infty $ (see Lemma VII.12.1
of \cite{Gamelin}). We denote by $\mathcal{C}$ the class of all cocycles on $%
K$ (identified modulo equality $\left( \mathfrak{m}_{\mathbb{R}}\,\times \,%
\mathfrak{m}_{K}\right) $-a.e. on $\mathbb{R}$\thinspace $\times \,K$), and
we associate with each $A\in \mathcal{C}$ the following one-parameter group $%
\mathcal{U}^{\left( A\right) }$ $\equiv \left\{ U_{t}^{\left( A\right)
}:t\in \mathbb{R}\right\} $ of linear bijections of $\mathcal{A}\left( 
\mathfrak{m}_{K}\right) $: for each $t\in \mathbb{R}$ and each $f\in 
\mathcal{A}\left( \mathfrak{m}_{K}\right) $, 
\begin{equation*}
\left( U_{t}^{\left( A\right) }f\right) \left( x\right) =A\left( t,x\right)
\,f\left( x+\mathfrak{e}_{t}\right) \text{, for }\mathfrak{m}_{K}\text{%
-almost all }x\in K\text{.}
\end{equation*}
It is readily seen that Theorem \ref{one-parameter-ver} applies to $\mathcal{%
U}^{\left( A\right) }$ $\equiv \left\{ U_{t}^{\left( A\right) }:t\in \mathbb{%
R}\right\} $. Briefly put, Helson's classic theory of generalized
analyticity and invariant subspaces uses the spectral decomposability of the
one-parameter unitary groups $\left\{ U_{t}^{\left( A\right) }\left|
L^{2}\left( \mathfrak{m}_{K}\right) \right. :t\in \mathbb{R}\right\} $, $%
A\in \mathcal{C}$, to establish a one-to-one correspondence between the
cocycles $A$ on $K$ and the normalized simply invariant subspaces of $%
L^{2}\left( \mathfrak{m}_{K}\right) $. This state of affairs has been
generalized to $L^{p}\left( \mathfrak{m}_{K}\right) $, $1\leq p<\infty $, as
follows (see \S \S\ 2,3 of \cite{ABG}): the one-parameter group $\mathcal{U}%
^{\left( A\right) }$ $\equiv \left\{ U_{t}^{\left( A\right) }:t\in \mathbb{R}%
\right\} $ transfers the Hilbert transform for $\mathbb{R}$ to a weak type $%
\left( 1,1\right) $ operator $\mathcal{H}^{\left( A\right) }:$ $L^{1}\left( 
\mathfrak{m}_{K}\right) \rightarrow \mathcal{A}\left( \mathfrak{m}%
_{K}\right) $ which is specified by taking%
\begin{equation*}
\left( \mathcal{H}^{\left( A\right) }f\right) \left( x\right)
=\lim_{n\rightarrow \infty }\int_{n^{-1}\leq \left| t\right| \leq n}\,\,%
\frac{\left( U_{-t}^{\left( A\right) }f\right) \left( x\right) }{\pi t}\,dt%
\text{, for }\mathfrak{m}_{K}\text{-almost all }x\in K\text{,}
\end{equation*}%
and $\mathcal{H}^{\left( A\right) }$ furnishes, via generalized Hardy
spaces, a cocycle characterization of the normalized simply invariant
subspaces of $L^{p}\left( \mathfrak{m}_{K}\right) $. For $1<p<\infty $, $%
\mathcal{H}^{\left( A\right) }\left| L^{p}\left( \mathfrak{m}_{K}\right)
\right. $ is a continuous linear mapping of $L^{p}\left( \mathfrak{m}%
_{K}\right) $ into itself. Clearly, Theorem \ref{one-parameter-ver} above,
when specialized to $\left\{ U_{t}^{\left( A\right) }:t\in \mathbb{R}%
\right\} $, provides the bilinear counterpart for $\mathcal{H}^{\left(
A\right) }$ (and, in contrast, the transferred bilinear Hilbert transform $%
\mathbb{H}_{\mathcal{U}^{\left( A\right) }}\left( \cdot ,\cdot \right) $ is
bounded even in the instances where the index $p_{3}$ of the target space $%
L^{p_{3}}\left( \mathfrak{m}_{K}\right) $ satisfies $\dfrac{2}{3}<p\leq $ $1$%
).
\end{example}

\label{aa_end_last_sec}


\begin{thebibliography}{99}
\bibitem{ABG} N. Asmar, E. Berkson, and T.A. Gillespie, \emph{Invariant
subspaces and harmonic conjugation on compact abelian groups, }Pacific J.
Math., \textbf{155}(1992), 201-213.

\bibitem{BBCG} E. Berkson, O. Blasco, M. Carro, and T.A. Gillespie, \emph{%
Discretization and transference of bisublinear maximal operators, }J.
Fourier Analysis and Appl., \textbf{12}(2006), 447-481.

\bibitem{Bou} J. Bourgain, \emph{Double recurrence and almost sure
convergence, }J. reine angew. Math., \textbf{404}(1990), 140-161.

\bibitem{Dem} C. Demeter, \emph{Pointwise convergence of the ergodic
bilinear Hilbert transform, }Ill J. Math, to appear.

\bibitem{DTT} C. Demeter, T. Tao, and C. Thiele, \emph{Maximal multilinear
operators,} Trans. Amer. Math. Soc., to appear.\emph{\ }

\bibitem{Doob} J.L. Doob, \emph{Stochastic Processes, }Wiley and Sons, New
York, 1953.

\bibitem{Doob1} J.L. Doob, \emph{A ratio operator limit theorem,} Z.
Wahrscheinlichkeitstheorie und Verw. Gebiete, \textbf{1(}1962/1963),
288--294.

\bibitem{Frem} D.H. Fremlin, \emph{Measure Theory. Vol. 3: Measure Algebras,}
T. Fremlin publ., Colchester, England, 2004.

\bibitem{Gamelin} T.W. Gamelin, \emph{Uniform Algebras,} Prentice-Hall,
Englewood Cliffs, N.J., 1969.

\bibitem{Halmos} P.R. Halmos, \emph{Measure Theory, }Van Nostrand,
Princeton, N.J., 1950.

\bibitem{PRH-vN} P.R. Halmos and J. von Neumann, \emph{Operator Methods in
Classical Mechanics, II,} Annals of Math, 2nd Ser., \textbf{43}(1942),
332-350.

\bibitem{Helson} H. Helson, \emph{Analyticity on compact abelian groups,} in
Algebras in Analysis, Proceedings of 1973 Birmingham Conference-NATO
Advanced Study Institute, Academic Press, London, 1975, pp. 1-62.

\bibitem{Lacey} M. Lacey, \emph{The bilinear maximal functions map into }$%
L^{p}$ \emph{for }$\dfrac{2}{3}<p\leq 1$\emph{, }Annals of Math., \textbf{151%
}(2000), 35-57.

\bibitem{Lamp} J. Lamperti, \emph{On the isometries of certain function}-%
\emph{spaces, }Pacific J. Math., \textbf{8}(1958), 459-466.

\bibitem{Rota} G-C Rota, \emph{An ``Alternierende Verfahren'' for general
positive operators, }Bull. Amer. Math. Soc., \textbf{68}(1962), 95-102.

\bibitem{Stein} E.M. Stein, \emph{Topics in Harmonic Analysis Related to the
Littlewood-Paley Theory, }Annals of Math. Studies 63, Princeton Univ. Press,
1970.
\end{thebibliography}
\end{document}